\documentclass[10pt]{article}
\usepackage{amsmath, amsthm,  latexsym, amssymb, amsfonts, epsfig,
caption2, graphicx}
\usepackage[dvipdfm]{hyperref}
\usepackage{mathrsfs}

\allowdisplaybreaks

\def\cajita{\rule{5pt}{5pt}}

\renewcommand{\d}{\mathrm{d}}
\newcommand{\ep}{\varepsilon}
\newcommand{\ph}{\varphi}
\newcommand{\wh}{\widehat}

\newcommand{\wt}{\widetilde}
\newcommand{\pr}{\prime}

\newcommand{\R}{\mathbb{R}}

\newcommand{\E}{\mathbb{E}}
\newcommand{\PP}{\mathbb{P}}

\newcommand{\F}{\mathcal{F}}
\newcommand{\lt}{\left}
\newcommand{\rt}{\right}
\newtheorem{theorem}{Theorem}[section]

\newtheorem{lemma}[theorem]{Lemma}
\par

\par
\newtheorem{corollary}[theorem]{Corollary}
\par
\par
\newtheorem{remark}[theorem]{Remark}
\par

\par

\par

\begin{document}

\title{Regularity Properties of Viscosity Solutions of
Integro-Partial Differential Equations of
Hamilton-Jacobi-Bellman Type}
\author{ Shuai Jing\thanks{This work is
supported by European Marie Curie Initial Training Network (ITN)
project: $"$Deterministic and Stochastic Controlled Systems and
Application$"$, FP7-PEOPLE-2007-1-1-ITN, No. 213841-2, the National
Basic Research  Program of China (973 Program) grant No.
2007CB814900 (Financial Risk) and the NSF of China (No.11071144).}\\
\small{ D\'epartement de Math\'ematiques, Universit\'e de Bretagne
Occidentale,} \small{29285 Brest C\'edex, France} \\
\small{School of Mathematics, Shandong University, 250100, Jinan,
China}\\\small{E-mail:
\href{mailto:shuai.jing@univ-brest.fr}{shuai.jing@univ-brest.fr}}}
\maketitle
\begin{abstract}{% We study the regularity properties of
%a certain class of integro-partial differential equations of
%Hamilton-Jocobi-Bellman type with terminal condition,
%which can be interpreted through a stochastic control system,
%composed of a forward and a backward stochastic differential equation,
%both driven by  a Brownian motion and a compensated Poisson random measure.
%More precisely, we prove that, under appropriate assumptions,
%the viscosity solution of such equations is  jointly Lipschitz and
%jointly semiconcave in $(t,x)\in\Delta\times\R^d$, for all compact
%time intervals $\Delta$ excluding the terminal time. Our approach is based on
%the method of  time change for the Brownian motion and on Kulik's
%transformation for the Poisson random measure. It extends earlier
%works by Buckdahn, Cannarsa and Quincampoix (2010) and by Buckdahn,
%Huang and Li (2011) on regularity properties for viscosity solutions of
%Hamilton-Jacobi-Bellman partial differential equations without and with obstacle.
We study the regularity properties of integro-partial differential equations of Hamilton-Jocobi-Bellman type with terminal condition, which can be interpreted through a stochastic control system, composed of a forward and a backward stochastic differential equation, both driven by a Brownian motion and a compensated Poisson random measure. More precisely, we prove that, under appropriate assumptions, the viscosity solution of such equations is  jointly Lipschitz and jointly semiconcave in $(t,x)\in\Delta\times\R^d$, for all compact time intervals $\Delta$ excluding the terminal time. Our approach is based on the time change for the Brownian motion and on Kulik's transformation for the Poisson random measure. }
\end{abstract}
\bigskip

\bigskip


\section{Introduction}
%\newpage
\setcounter{equation}{0}

We are interested in the regularity properties of the viscosity solution
for a certain class of integro-partial differential equations (IPDEs) of
Hamilton-Jacobi-Bellman (HJB) type. In order to be more precise, let us consider the
following possibly degenerate equation:
\begin{equation}\label{3ipde}
\left\{
\begin{aligned}
\frac{\partial}{\partial t} V(t,x)+\inf_{u\in
U}\{(\mathcal{L}^u+B^u)V(t,x)+f\big(t,x,&V(t,x),(D_x V\sigma)(t,x), \\
& V(t, x+\beta(t,x,u,\cdot))-V(t, x),u\big)\}=0;\\
%\ \ (t,x)\in(0,T)\times\R^d;\\
V(T,x)=\Phi(x),\qquad \qquad\qquad\qquad\qquad\qquad &
%\\qquad\qquad
%x\in\R^d,
\end{aligned}
\right.
\end{equation}
where $U$ is a compact metric space,
$\mathcal{L}^u$ is the linear second order differential
operator
\[
\mathcal{L}^u
\ph(x)={\rm{tr}}\lt(\frac12\sigma\sigma^T(t,x,u)D_{xx}^2\ph(x)\rt)+b(t,x,u)\cdot
D_x \ph(x),\ \ \ \ph\in C^2(\R^d),
\]
and $B^u$ is the integro-differential operator:
\[
B^u \ph(x)=\int_E[\ph(x+\beta(t,x,u,e))-\ph(x)-\beta(t,x,u,e)\cdot
D_x\ph(x)]\Pi(\d e),\ \ \ \ph\in C_b^1(\R^d).
\]
Here, $\Pi$ denotes a finite L\'evy measure on $E=\R^n\backslash \{0\}.$
Our main results say that, under appropriate assumptions, for all $\delta>0$,
the viscosity solution $V$ is jointly Lipschitz and jointly semiconcave
on $[0,T-\delta]\times\R^d$, i.e., there is some constant $C_\delta$ such that
$$
|V(t_0,x_0)-V(t_1,x_1)|\le C_\delta(|t_0-t_1|+|x_0-x_1|),
$$
$$
\lambda V(t_0,x_0)+(1-\lambda)V(t_1,x_1)\le V\big(\lambda (t_0,x_0)+(1-\lambda)(t_1,x_1)\big)
+C_\delta\lambda(1-\lambda)(|t_0-t_1|^2+|x_0-x_1|^2),
$$
for all $(t_0,x_0), (t_1,x_1)\in[0,T-\delta]\times\R^d$.
The joint semiconcavity of $V$ stems its importance from the fact that,
due to Alexandrov's theorem,  it implies that $V$ has a
second order expansion in $(t,x)$, $\d t\d x$-a.e.
Such expansions are important, for instance,
for the study of the propagation of singularities.

Although, at least for PDEs of HJB type, the regularity of
the solution of strictly elliptic equations
(with $\sigma\sigma^T\ge \alpha I$, for $\alpha>0$)
has been well understood  for a long time,  the joint
regularity (Lipschitz continuity and semiconcavity) in
$(t,x)$ for the viscosity solution of such equations,
for which $\sigma\sigma^T$ is not necessarily strictly
elliptic, have been studied only recently.
However, under suitable hypotheses,
the Lipschitz continuity and the semiconcavity  of $V(t,x)$ in $x$
as well as the H\"older continuity of $V(t,x)$ in $t$
has already been known for a longer time.
As concerns the semiconcavity of $V(t,x)$ in $x$,
a purely analytical proof was given by Ishii and Lions \cite{IL};
for a stochastic proof the reader is referred, for example,
to Yong and Zhou \cite{YoZh}.
Concerning the Lipschitz continuity of $V$ in $x$ and
the H\"older continuity in $t$
(with H\"older coefficient $1/2$),
we refer, for example, to Pham \cite{Ph}.
Krylov \cite{Kr} suggested the joint Lipschitz continuity of
$V(t,x)$ on $(t,x)$. However, counterexamples show that,
in general, one cannot get the Lipschitz continuity or
semiconcavity in $(t,x)$ for the whole domain $[0,T]\times\R^d$.
In Buckdahn, Cannarsa and Quincampoix \cite{BCQ},
it was shown that the viscosity solution of PDEs of HJB type
(with $\beta=0$) is Lipschitz and semiconcave over
$[0,T-\delta]\times\R^d$ for $\delta>0$.
In Buckdahn, Huang and Li \cite{BHL},
these results were extended to PDEs with obstacle.
The approach in \cite{BCQ} and \cite{BHL} consists
in the study of the viscosity solution $V$ with the
help of its stochastic interpretation as a value
function of an associated stochastic control problem;
it uses, in particular, the method of time change,
which was translated to backward stochastic differential equations (BSDEs) in \cite{BHL}.

In this paper we study the joint regularity of $V(t,x)$ in $(t,x)$
through the stochastic interpretation of the above HJB equation
as a stochastic control problem composed of a forward and
a backward stochastic differential equation (SDE).
More precisely, let $(t,x)\in[0,T]\times\R^d$,
$B=(B_s)_{s\in[t,T]}$ be a $d$-dimensional Brownian motion
with initial value zero at time $t$, and let $\mu$ be a
Poisson random measure on $[t,T]\times E$.
We denote by $\mathbb{F}$ the filtration generated by $B$ and $\mu$,
and by $\mathcal{U}^{B,\mu}(t,T)$ the set of all
$\mathbb{F}$-predictable control processes with values in $U$.
It is by now standard that the SDE
driven by the Brownian motion $B$ and the compensated Poisson random measure
$\tilde{\mu}$:
\begin{equation}\label{3eq:sdej}
X_s^{t,x,u}=x+\int^s_t\hspace{-1mm}b(r,X_r^{t,x,u},u_r)\d
r+\int^s_t\hspace{-1mm}\sigma(r,X_r^{t,x,u},u_r)\d
B_r+\int^s_t\hspace{-2mm}\int_E\beta(r,X_{r-}^{t,x,u},u_r,e)\tilde{\mu}(\d r, \d e), \ s\in[t,T],
\end{equation}
has a unique solution under appropriate assumptions for the coefficients.
With this SDE we associate the BSDE with jumps
\begin{eqnarray}\label{3bsde}
Y_s^{t,x,u}&=&\Phi(X_T^{t,x,u}) +\int^T_s
f(r,X_r^{t,x,u},Y_r^{t,x,u},Z_r^{t,x,u},U_r^{t,x,u},u_r)\d
r -\int^T_s Z_r^{t,x,u}\d B_r\nonumber\\
&&-\int^T_s\int_EU_r^{t,x,u}(e)\tilde{\mu}(\d r,
\d e), \ \ s\in[t,T].
\end{eqnarray}
(As concerns the assumptions on the coefficients, we refer to the hypotheses {\bf(H1)-(H5)}
in Section 2 and Section 3.)
From Barles, Buckdahn and Pardoux \cite{BBP} we know that the above
BSDE with jumps (\ref{3bsde}) has a unique square integrable
solution $(Y^{t,x,u},Z^{t,x,u},U^{t,x,u})$.
Moreover, since $Y^{t,x,u}$ is $\mathbb{F}$-adapted,
$Y^{t,x,u}_t$ is deterministic. It follows from Barles, Buckdahn and Pardoux  \cite{BBP}
or Pham \cite{Ph} that the value function
\begin{equation}\label{eq:value}
V(t,x)=\inf_{u\in\mathcal{U}^{B,\mu}(t,T)}Y^{t,x,u}_t, \ \ (t,x)\in[0,T]\times\R^d
\end{equation}
is the viscosity solution of our IPDE.

Since unlike \cite{BCQ} and \cite{BHL}, our system involves not only
the Brownian motion $B$ but also the Poisson random measure $\mu$,
the method of time change for the Brownian motion alone is not sufficient for
our approach here. So we combine the method of time change for the
Brownian motion by Kulik's transformation
for Poisson random measures (see, \cite{Kul0} \cite{Kul}).
To our best knowledge, the use of Kulik's transformation for the study
of stochastic control problems is new.
Because of the difficulty to obtain suitable $L^p$-estimates
of the stochastic integrals with respect the compensated
Poisson random measure (see, for example, Pham \cite{Ph})
we have to restrict ourselves to the case of a finite
L\'evy measure $\Pi(E)<+\infty$. The more general case
where $\int_E(1\wedge |e|^2)\Pi(\d e)<+\infty$ remains still open.

Our paper is organized as follows. In Section 2 we introduce our
main tools, i.e., the method of time change for the Brownian motion
and Kulik's transformation for the Poisson random measure,
with the help of which we study the joint Lipschitz continuity
for the viscosity solution of the IPDEs of HJB type.
This method of time change for the Brownian motion combined
with Kulik's transformation is extended in Section 3 to
the study of the semiconcavity property for the viscosity solution of
IPDE (\ref{3ipde}). The proof of more technical statements and
estimates used in Section 3 is shifted in the Appendix.

\section{Lipschitz Continuity}
\setcounter{equation}{0}
In this section, we prove the joint Lipschitz
continuity of the viscosity solution of a certain
class of integro-differential Hamilton-Jacobi-Bellman
(HJB) equations.

Let $T$ be an arbitrarily fixed time horizon,
$U$ a compact metric space, $E=\R^d\backslash\{0\}$
and $\mathcal{B}(E)$ be the Borel $\sigma$-algebra over $E$.
We are concerned with the integro-partial
differential equation of HJB type (\ref{3ipde}).
%\begin{equation}\label{eq:hjb}
%\left\{
%\begin{array}{lr}
%\frac{\partial}{\partial t} V(t,x)+\inf_{u\in
%U}\{(\mathcal{L}^u+B^u)V(t,x)
%+f(t,x,V(t,x),(D_x V\sigma)(t,x), \\
%\qquad\qquad\qquad\qquad V(t, x+\beta(t,x,u,\cdot))-V(t, x),u)\}=0,
%\ \
%(t,x)\in(0,T)\times\R^d;\\
%V(T,x)=\Phi(x),\qquad \qquad\qquad\qquad\qquad\qquad\qquad\quad\qquad\qquad\qquad\ \qquad x\in\R^d,
%\end{array}
%\right.
%\end{equation}
%where, for $(t,x,u)\in[0,T]\times\R^d\times U$,
%$\mathcal{L}^u$ is the linear second order differential
%operator
%\[
%\mathcal{L}^u
%\ph(x)={\rm{tr}}\lt(\frac12\sigma\sigma^T(t,x,u)D_{xx}^2\ph(x)\rt)+b(t,x,u)\cdot
%D_x \ph(x),\ \ \ph\in C^2(\R^d),
%\]
%$B^u$ is the integro-differential operator
%\[
%B^u \ph(x)=\int_E[\ph(x+\beta(t,x,u,e))-\ph(x)-\beta(t,x,u,e)\cdot
%D_x\ph(x)]\Pi(\d e), \ \ \ph\in C^1_b(\R^d),
%\]
%and $\Pi$ is a finite L\'evy measure defined on the measurable space $(E,\mathcal{B}(E))$.
The coefficients
\[
b:[0,T]\times \R^d\times U\to \R^d,\ \sigma:
[0,T]\times \R^d\times U\to \R^{d\times d},\
\beta:[0,T]\times\R^d\times U\times E\to\R^d,
\]
\[
f:[0,T]\times\R^d\times\R\times\R^d\times L^2(E,\mathcal{B}(E),\Pi;\R)\times U\to\R
\ \ \rm{and}\ \ \Phi:\R^d\to \R
\]
are bounded continuous functions which satisfy the following conditions:

\noindent(\textbf{H1}) There
exists a constant $K>0$ such that, for any
$\xi_i=(s_i,y_i)\in[0,T]\times\R^d, u\in U, i=1,2$,
\[
\begin{aligned}
&|b(\xi_1,u)-b(\xi_2,u)|+|\sigma(\xi_1,u)-\sigma(\xi_2,u)|
+\lt(\int_E|\beta(\xi_1,u,e)-\beta(\xi_2,u,e)|^4\Pi(\d e)\rt)^{1/4}\\
\le&
K(|s_1-s_2|+|y_1-y_2|).
\end{aligned}
\]

\noindent(\textbf{H2}) The function $f$ is Lipschitz in $(t,x,y,z,p)$,
uniformly with respect to $u\in U$, and the function
$\Phi$ is a Lipschitz function.

The integro-PDE (\ref{3ipde}), as is well-known by now
(see, for instance, \cite{BBP}), has a unique continuous viscosity solution $V(t,x)$
in the class of the continuous functions with at most polynomial growth.

Let $\{B_s^0\}_{s\ge0}$ be a $d$-dimensional Brownian motion
defined on a complete space $(\Omega_1,\F_1,$ $\PP_1)$,
and $\eta$ be a Poisson random measure defined
on a complete probability space  $(\Omega_2, \F_2,$ $\PP_2)$.
We introduce $(\Omega, \F, \PP)$ as the product space $(\Omega,\F,\PP) =
(\Omega_1,\F_1,\PP_1) \otimes(\Omega_2,\F_2,\PP_2)
=(\Omega_1\times\Omega_2,\F_1\otimes\F_2,\PP_1\otimes\PP_2)$.
The processes $B^0$ and $\eta$ are canonically extended from
$(\Omega_1,\F_1,$ $\PP_1)$ and $(\Omega_2, \F_2,$ $\PP_2)$, respectively,
to the product space $(\Omega, \F, \PP)$.
We denote the compensated Poisson random measure associated with $\eta$ by
$\tilde{\eta}$, i.e., $\tilde{\eta}(\d t, \d e)=\eta(\d t, \d e)-\d
t\Pi(\d e)$. We assume throughout this paper that the L\'evy measure
$\Pi$ is a finite measure on $(E,\mathcal{B}(E))$.

We define the process $\{B_s\}_{s\ge t}$ by putting
\begin{equation}\label{BM}
B_s=B_s^0-B^0_{t}, \ \ s\in[t,T],
\end{equation} so that $\{B_s\}_{s\ge t}$ is a Brownian motion
beginning at time $t$ with $B_t=0$. Furthermore, we denote by $\mu$
be the restriction of the Poisson random measure $\eta$
from $[0,T]\times E$ to $[t,T]\times E$,
and by $\tilde{\mu}$ its compensated measure.

We put
$$
\mathcal{F}_s^B=\sigma\{B_r, r\in[t,s]\}\vee
\mathcal{N}_{\mathbb{P}_1}, \ \
\F_s^\mu=\sigma\{\mu((t,r]\times\Delta):
\Delta\in\mathcal{B}(E),r\in[t,s]\}\vee \mathcal{N}_{\mathbb{P}_2},
$$
and
$$
\F_s=\lt(F_s^B\otimes\F_s^\mu\rt)\vee\mathcal{N}_{\PP},\ \ s\in[t,T],
$$
where $\mathcal{N}_{\mathbb{P}_1}$, $\mathcal{N}_{\mathbb{P}_2}$
and $\mathcal{N}_{\mathbb{P}}$ are the collections of the null sets
under the corresponding probability measure.

Let us also introduce the following spaces of stochastic processes
over $(\Omega,\F,\PP)$ which will be needed in what follows.
By $\mathcal{S}^2(t,T;\R^d)$ we denote the set of all
$\mathbb{F}$-adapted c\`adl\`ag processes $\{Y_s;t\le s\le T\}$ such that
$$\|Y\|_{\mathcal{S}^2(t,T;\R^d)}=\E\lt[\sup_{t\le s\le T}|Y_s|^2\rt]< \infty.$$
Let $\mathbb{L}^2(t,T;\R^d)$ denote the set of all
$\mathbb{F}$-predictable $d$-dimensional processes $\{Z_s: t\le s\le T\}$ such that
\[
\|Z\|_{\mathbb{L}^2(t,T;\R^d)}=\lt(\E\lt[\int^T_t|Z_s|^2\d s\rt]\rt)^{1/2}<\infty.
\]
Finally, we also introduce the space $\mathbb{L}^2(t,T;\tilde{\mu},\R)$
of mappings $U:\Omega\times[0,T]\times E\to\R$ which are
$\mathbb{F}$-predictable and measurable such that
\[
\|U\|_{\mathbb{L}^2(t,T;\tilde{\mu},\R)}=
\lt(\E\lt[\int^T_t\int_E |U_s(e)|^2\Pi(\d e)\d s\rt]\rt)^{1/2}<\infty.
\]

Let us now consider the following stochastic differential equation driven by
the Brownian motion $B$ and the compensated Poisson random measure
$\tilde{\mu}$:
\begin{equation}\label{eq:sdej}
X_s^{t,x,u}=x+\int^s_t\hspace{-1mm}b(r,X_r^{t,x,u},u_r)\d
r+\int^s_t\hspace{-1mm}\sigma(r,X_r^{t,x,u},u_r)\d
B_r+\int^s_t\hspace{-2mm}\int_E\beta(r,X_{r-}^{t,x,u},u_r,e)\tilde{\mu}(\d r, \d e),\ s\in[t,T],
\end{equation}
where the process $u:[t,T]\times\Omega\to U$ is an admissible
control, i.e., an $\mathbb{F}$-predictable process with values in $U$;
the space of admissible controls over the time interval $[t,T]$
is denoted by  $\mathcal{U}^{B,\mu}(t,T)$. The following theorem is by now
classical:
\begin{theorem}
Assume the Lipschitz condition $(\textbf{H1})$. For any fixed admissible
control $u(\cdot)\in\mathcal{U}(t,T)$, there exists a unique adapted
c\`adl\`ag solution $(X_s^{t,x,u})_{s\in[t,T]}\in\mathcal{S}^2(t,T;\R^d)$ of the stochastic
differential equation $(\ref{eq:sdej})$.
\end{theorem}

We associate SDE (\ref{eq:sdej}) with the backward stochastic differential equation
\begin{eqnarray*}
Y_s^{t,x,u}&=&\Phi(X_T^{t,x,u}) +\int^T_s
f(r,X_r^{t,x,u},Y_r^{t,x,u},Z_r^{t,x,u},U_r^{t,x,u},u_r)\d
r -\int^T_s Z_r^{t,x,u}\d B_r\nonumber\\
&&-\int^T_s\int_EU_r^{t,x,u}(e)\tilde{\mu}(\d r,
\d e), \ \ s\in[t,T].
\end{eqnarray*}
Then from Barles, Buckdahn and Pardoux \cite{BBP}, Tang and Li \cite{TL},
we know that this BSDE has a unique solution
\[
(Y^{t,x,u},Z^{t,x,u},U^{t,x,u})\in \mathcal{S}^2(t,T;\R)
\times\mathbb{L}^2(t,T;\R^d) \times\mathbb{L}^2(t,T;\tilde{\mu},\R).
\]
Notice that $Y^{t,x,u}_t$ is $\F_t$-measurable,
hence it is deterministic in the sense that it coincides $\PP$-a.s.
with a real constant, with which it is identified.  Thus, we have
$$
Y_t^{t,x,u}=\E\lt[Y_t^{t,x,u}\rt]=
\E\lt[\int^T_tf(r,X_r^{t,x,u},Y_r^{t,x,u},Z_r^{t,x,u},U_r^{t,x,u},u_r)\d
r+\Phi(X_T^{t,x,u})\rt].
$$
As usual in stochastic control problems, we define
the cost functional $J(t,x;u)$ associated with
$u\in\mathcal{U}^{B,\mu}(0,T)$ by setting $J(t,x;u):=Y^{t,x,u}_t$,
and the value function is defined as follows:
\[
V(t,x)=\inf_{u(\cdot)\in\mathcal{U}^{B,\mu}(t,T)}J(t,x;u), \ \
(t,x)\in[0,T]\times\R^d.
\]
It is well known by now that  $V=\{V(t,x): (t,x)\in[0,T]\times\R^d\}$ is a continuous
viscosity solution of the HJB equation (\ref{3ipde}).
Moreover, $V$ is the unique viscosity solution in the
class of continuous functions with at most polynomial growth
(see: \cite{Ph}, \cite{Wu}).

Our main result in this section is the following theorem.
\begin{theorem}\label{thm:lip}
Let $\delta\in (0,T)$ be arbitrary but fixed.
Under our assumptions $(\textbf{H1})$ and $(\textbf{H2})$,
the value function
$V(\cdot,\cdot)$ is jointly Lipschitz continuous on
$[0,T-\delta]\times\R^d$,  i.e., for some constant $C_\delta$ we have,
for all $(t_0,x_0), (t_1,x_1)\in[0,T-\delta]\times\R^d$:
$$
|V(t_0,x_0)-V(t_1,x_1)|\le C_\delta(|t_0-t_1|+|x_0-x_1|).
$$
\end{theorem}
\begin{remark}\label{remark_lip}
In general we cannot expect to get the joint Lipschitz continuity over the whole
domain $[0,T]\times\R^d$. In \cite{BCQ} is given an easy counterexample: We study the problem
\[
X_s^{t,x}=x+B_s, \ \ s\in[t,T],\  \ x\in\R;
\]
\[
Y^{t,x}_s=-\E\lt[\lt|X_T^{t,x}\rt||\F_s\rt]=-\E\lt[|x+B_T||\F_s\rt],
\ \ s\in[t,T],
\]
without control neither jumps. Then
\[
V(t,x)=Y^{t,x}_t=-\E\lt[|x+B_T|\rt],
\]
and, for $x=0$, recalling that $B$ is a Brownian motion with $B_t=0$, we have
\[
V(t,0)=-\E[|B_T|]=-\sqrt{\frac{2}{\pi}}\sqrt{T-t},\ \ t\in[0,T].
\]
Obviously, $V(\cdot,x)$ is not Lipschitz in $t$ for $t=T$. However,
$V$ is jointly Lipschitz on $[0,T-\delta]\times\R,$ for $\delta\in(0,T)$.
\end{remark}

Let us introduce now Kulik's transformation in our framework. The reader interested
in more details on this transformation is referred to the papers \cite{Kul0} and \cite{Kul}.

Let $t_0,t_1\in[0,T]$ and  let, for $t=t_0$, $\mu$ be the Poisson random measure
which we have introduced as restriction of $\eta$ from $[0,T]\times E$ to $[t_0,T]\times E$.
With the help of $\mu$ we define now a random measure $\tau(\mu)$ on $[t_0,T]\times E$. Denoting by
\[
\tau: [t_1,T]\to[t_0,T]
\]
the linear time change
\[
\tau(s)=t_0+\frac{T-t_0}{T-t_1}(s-t_1), \ s\in[t_1,T],
\]
we put
\[
\tau(\mu)([t_1,s]\times \Delta):=\mu([\tau(t_1),\tau(s)]\times\Delta),
\ t_1\le s\le T,\ \Delta\in\mathcal{B}(E).
\]
Observing that $\dot{\tau}=\dot{\tau}(s)=\frac{T-t_0}{T-t_1}$,
we put
\[\gamma=\ln \lt({\dot{\tau}}\rt)=\ln \lt(\frac{T-t_0}{T-t_1}\rt).\]
From Lemma 1.1 in Kulik \cite{Kul0} we know that, for all $\{s_1, \cdots,
s_n\}\subset[t_1,T]$, $\Delta_1,\cdots,\Delta_n\in\mathcal{B}(E)$ and
all Borel function $\varphi:\R^n\to\R^+$, we have
\begin{lemma}
\[
\begin{aligned}
&\E\lt[\ph(\tau(\mu)([t_1,s_1]\times\Delta_1),\cdots, \tau(\mu)([t_1,s_n]\times\Delta_n))\rt]\\
=&\E\lt[\rho_{\tau}\ph(\eta([t_1,s_1]\times\Delta_1),\cdots, \eta([t_1,s_n]\times\Delta_n))\rt],
\end{aligned}
\]
where
\[
\rho_\tau=\exp\{\gamma\eta([t_1,T]\times
E)-(t_1-t_0)\Pi(E)\}.
\]
\end{lemma}
For the convenience of the reader we sketch the proof.
However, we restrict to a special case ($n=1$),
the proof of the general case $n\ge 1$ can be carried out with a similar argument
and can be consulted for the more general case $\Pi(E)=+\infty$ in \cite{Kul0}.
\bigskip

\noindent\textit{Proof}: (for $n=1$). Observing that for $\Delta_2=E\backslash \Delta_1$,
\[
\eta([t_1,s_1]\times\Delta_1), \ \ \ \ \eta([s_1,T]\times\Delta_1)
\ \ \ \ \mbox{and}\ \ \ \ \eta([t_1,T]\times\Delta_2)
\]
are independent Poisson distributed random variables with the
intensities $(s_1-t_1)\Pi(\Delta_1)$, $(T-s_1)\Pi(\Delta_1)$
and $(T-t_1)\Pi(\Delta_2)$, respectively,
we have
\[
\begin{aligned}
&\E[\rho_{\tau}\ph(\eta([t_1,s_1]\times\Delta_1))]\\
=&\bigg\{\sum_{k,l\ge0}\ph(k)\exp\{\gamma k+\gamma l-(t_1-t_0)
\Pi(\Delta_1)\}\times\\
&\qquad\qquad\times\exp\{-(T-t_1)\Pi(\Delta_1)\}
\frac{((s_1-t_1)\Pi(\Delta_1))^k}{k!} \frac{((T-s_1)\Pi(\Delta_1))^l}{l!}\bigg\} \\
&\times\bigg\{\sum_{m\ge0}\exp\{\gamma m-(t_1-t_0)\Pi(\Delta_2)\}\exp\{-(T-t_1)
\Pi(\Delta_2)\}\frac{((T-t_1)\Pi(\Delta_2))^m}{m!}\bigg\}\\
=&I_1\times I_2.
\end{aligned}
\]
But, taking into account the definition of $\gamma$ and that of $\tau(s_1)$ we have
\[
\begin{aligned}
I_1=&\bigg\{\sum_{k\ge0}\ph(k)\exp\{-(T-t_0)\Pi(\Delta_1)\}\frac{1}{k!}
\lt(\lt(\frac{T-t_0}{T-t_1}\rt)(s_1-t_1)\Pi(\Delta_1)\rt)^k\bigg\}\times\\
&\qquad\qquad\times\bigg\{\sum_{l\ge0}\frac{1}{l!}
\lt(\lt(\frac{T-t_0}{T-t_1}\rt)(T-s_1)\Pi(\Delta_1)\rt)^l\bigg\}
\\
=&\bigg\{\sum_{k\ge0}\ph(k)\exp\{-(T-t_0)\Pi(\Delta_1)\}\frac{1}{k!}
\lt((\tau(s_1)-t_0)\Pi(\Delta_1)\rt)^k\bigg\}\exp\{(T-\tau(s_1))\Pi(\Delta_1)\}\\
=&\sum_{k\ge0}\ph(k)\exp\{-(\tau(s_1)-t_0)\Pi(\Delta_1)\}
\frac{1}{k!}\lt((\tau(s_1)-t_0)\Pi(\Delta_1)\rt)^k\\
=&\E[\ph(\tau(\mu)([t_1,s_1]\times\Delta_1))].
\end{aligned}
\]
In analogy to the computation for $I_1$, but now with $\ph\equiv1$,
we get that $I_2=1$. Consequently,
\[
\E[\rho_{\tau}\ph(\eta([t_1,s_1]\times\Delta_1))]=
\E[\ph(\tau(\mu)([t_1,s_1]\times\Delta_1))].
\]
Hence the proof is complete.
\hfill$\cajita$
\bigskip

From the above lemma we have, for all $n\ge1$, $\{s_1,\cdots,s_n\}\subset[t_1,T],
\Delta_1,\cdots, \Delta_n\in\mathcal{B}(E)$ and $\ph:\R^n\to\R^+$ Borel function,
\[
\begin{aligned}
&\E[\ph(\eta([t_1,s_1]\times\Delta_1),\cdots,\eta([t_1,s_n]\times\Delta_n))]\\
=&\E\big[\rho_{\tau}\ph(\eta([t_1,s_1]\times\Delta_1),\cdots,\eta([t_1,s_n]\times\Delta_n))\\
&\qquad\times\exp\{-\gamma\eta([t_1,T]\times E)+(t_1-t_0)\Pi(E)\}\big]\\
=&\E\big[\ph(\tau(\mu)([t_1,s_1]\times\Delta_1),\cdots,\eta([t_1,s_n]\times\Delta_n))\\
&\qquad\times\exp\{-\gamma\tau(\mu)([t_1,T]\times E)+(t_1-t_0)\Pi(E)\}\big]\\
=&\E\lt[g_{\tau}\ph(\tau(\mu)([t_1,s_1]\times\Delta_1),\cdots,\eta([t_1,s_n]\times\Delta_n))\rt],
\end{aligned}
\]
for $g_{\tau}=\exp\{-\gamma\tau(\mu)([t_1,T]\times E)+(t_1-t_0)\Pi(E)\}$.

This allows to show that under the probability measure $\mathbb{Q}_{\tau}=g_{\tau}\PP$,
the point process $\tau(\mu)$ defined over $[t_1,T]\times E$, has the same law as
the Poisson random measure $\eta$ restricted to $[t_1,T]\times E$, under $\PP$.
Consequently, under $\mathbb{Q}_{\tau}=g_{\tau}\PP$, $\tau(\eta)$ is a Poisson
random measure with compensator $\d s \Pi (\d e)$.

We use the same time change $\tau:[t_1,T]\to[t_0,T]$
in order to introduce the process
\[
W_t=\frac{1}{\sqrt{\dot{\tau}}}B_{\tau(t)},\ t\in[t_1,T].
\]
We observe that $W=(W_t)_{t\in[t_1,T]}$ is a Brownian motion under
the probability $\PP$ but also under $\mathbb{Q}_{\tau}=g_{\tau}\PP$
(Indeed, $B$ and $g_\tau$ are independent under $\PP$),
and $W$ and $\tau(\mu)$  are independent under both $\PP$ and $\mathbb{Q}_{\tau}$.

Let $\ep>0$. From the definition of the value function $V$,
\[
V(t_0,x_0)=\inf_{u\in\mathcal{U}^{B,\mu}(t_0,T)}J(t_0,x_0,u),
\]
we get the existence of an admissible control $u^0\in\mathcal{U}^{B,\mu}(t_0,T)$
such that
\[
J(t_0,x_0,u^0)\le V(t_0,x_0)+\ep.
\]
We define
\[
u^1(t)=u^0(\tau(t)),\ t\in[t_1,T].
\]
Then, obviously, $u^1\in\mathcal{U}^{W,\tau(\mu)}(t_1,T)$, i.e.,
$u^1$ is a $U$-valued process predictable with respect the filtration
\[
\F_{t}^{W,\tau(\mu)}=\sigma\{W_s,\tau(\mu)([t_1,s]\times\Delta), s\in[t_1,t],
\Delta\in\mathcal{B}(E)\}\vee \mathcal{N}_{\PP}, \ \ t\in[t_1,T],
\]
generated by $W$ and $\tau(\mu)$.

Let now
$X^0=\lt\{X_s^0\rt\}_{s\in[t_0,T]}$ be the solution of the forward equation
\begin{equation}\label{eq:sdej1}
X_s^0=x_0+\int^s_{t_0}b(r,X_r^0,u^0_r)\d
r+\int^s_{t_0}\sigma(r,X_r^0,u^0_r)\d
B_r+\int^s_{t_0}\int_E\beta(r,X_{r-}^0,u^0_r,e)\tilde{\mu}(\d
r, \d e),\ \ s\in[t_0,T],
\end{equation}
under the probability $\PP$,
and let $X^1=\lt\{X_s^{1}\rt\}_{s\in[t_1,T]}$ be the solution of the equation
\begin{equation}\label{eq:sdej2}
X_s^1=x_1+\int^s_{t_1}b(r,X_r^1,u^1_r)\d
r+\int^s_{t_1}\sigma(r,X_r^1,u^1_r)\d
W_r+\int^s_{t_1}\int_E\beta(r,X_{r-}^1,u^1_r,e)\widetilde{\tau(\mu)}^{\mathbb{Q}_{\tau}}(\d
r, \d e),\ \ s\in[t_1,T],
\end{equation}
under probability measure
$\mathbb{Q}_{\tau}$.
Notice that the compensated Poisson random measure
$\tilde{\mu}$ under $\PP$ is of the form
\[
\tilde{\mu}(\d s, \d e)=\mu(\d s, \d e)-\d s \Pi(\d e),
\ (s,e)\in[t_0,T]\times E,
\]
while the compensated Poisson random measure for $\tau(\mu)$
under $\mathbb{Q}_{\tau}$ has the form
\[
\wt{\tau(\mu)}^{\mathbb{Q}_{\tau}}(\d s, \d e)=
\tau(\mu)(\d s, \d e)-\d s \Pi(\d e),  \ (s,e)\in[t_1,T]\times E.
\]

We employ the BSDE method to prove the Lipschitz continuity of the value function $V$.
For this we associate the above SDEs with the following BSDEs with jumps:
\begin{equation}\label{eq:bsdej1}
Y^0_s=\Phi(X_T^0)+\int^T_{s}f(r,X_r^0,Y_r^0,Z_r^0,U_r^0,u^0_r)\d
r-\int^T_s Z_r^0\d B_r-\int^T_s\int_E
U^0_r(e)\tilde{\mu}(\d r,\d e), \ s\in[t_0,T],
\end{equation}
under probability $\PP$, and
\begin{equation}\label{eq:bsdej2}
Y^1_s=\Phi(X_T^1)+\int^T_{s}f(r,X_r^1,Y_r^1,Z_r^1,U_r^1,u^1_r)\d
r-\int^T_s Z_r^1\d W_r-\int^T_s\int_E
U^1_r(e)\widetilde{\tau(\mu)}^{\mathbb{Q}_{\tau}}(\d r,\d e),\ s\in[t_1,T],
\end{equation}
under probability $\mathbb{Q}_{\tau}$.
From \cite{BBP}, we know the above two
BSDEs have unique solutions $(Y^0,$
$Z^0,$ $ U^0)$$=(Y_s^0, Z_s^0, U_s^0)_{s\in[t_0,T]}$ and
$(Y^1,Z^1,U^1)=(Y_s^1, Z_s^1, U_s^1)_{s\in[t_1,T]}$, respectively.
While $Y^0$ is adapted and $Z^0$ and $U^0$ are predictable
with respect to the filtration generated by $B$ and $\mu$,
$Y^1$ is adapted and $Z^1$ and $U^1$ are predictable
with respect the filtration generated by $W$ and $\tau(\mu)$.
Thus, $Y^0_{t_0}$ and $Y^1_{t_1}$ are deterministic, and
from the definition of the cost functionals we have
\begin{equation}\label{cost1}
Y_{t_0}^{0}-\ep=J(t_0,x_0;u^0)-\ep\le V(t_0,x_0)
\lt(=\inf_{u\in\mathcal{U}^{B,\mu}(t_0,T)}J(t_0,x_0;u)\rt),
\end{equation}
and
\begin{equation}\label{cost2}
Y_{t_1}^{1}=J(t_1,x_1;u^1)\ge V(t_1,x_1)
\lt(=\inf_{u\in\mathcal{U}^{W,\tau(\mu)}(t_1,T)}J(t_1,x_1;u)\rt).
\end{equation}
Here we have used that the stochastic interpretation of $V$
does not depend on the special choice of the underlying
driving Brownian motion and the underlying Poisson random measure
with compensator $\d s \Pi(\d e)$.
In order to show the Lipschitz property of $V$ in $(t,x)$, we have to estimate
\[
\begin{aligned}
V(t_0,x_0)-V(t_1,x_1)\ge& J(t_0,x_0;u^0)-J(t_1,x_1;u^1)-\ep\\
=&Y_{t_0}^0-Y_{t_1}^1-\ep.
\end{aligned}
\]
However, in order to estimate the difference between the processes $Y^0$ and $Y^1$,
we have to make their both BSDEs comparable, i.e.,
we need them over the same time interval,
driven by the same Brownian motion and by the same compensated Poisson random measure.
For this reason we apply to SDE (\ref{eq:sdej2}) and BSDE (\ref{eq:bsdej2})
the inverse time change $\tau^{-1}: [t_0,T]\to[t_1,T]$. So we introduce the process
$\wt{X}^1=\{\wt{X}_s^1\}_{s\in[t_0,T]}$ by setting
$\wt{X}_s^1=X^1_{\tau^{-1}(s)}$.
We also observe that
$W_{\tau^{-1}(r)}=\frac{1}{\sqrt{\dot{\tau}}}B_r$
and
$u^1_{\tau^{-1}(r)}=u^0_r$, $r\in[t_0,T]$.
Obviously,
$\wt{X}^1\in \mathcal{S}^2(t_0,T;\R)$
is the unique solution of the SDE
\begin{eqnarray}\label{eq:sdej3}
\wt{X}_s^1
&=&x_1+\int^s_{t_0}
b(\tau^{-1}(r),\wt{X}_r^1,u^1_{\tau^{-1}(r)})\d
\tau^{-1}(r)+\int^s_{t_0}\sigma(\tau^{-1}(r),
\wt{X}_r^1,u^1_{\tau^{-1}(r)})\d
W_{\tau^{-1}(r)}\nonumber\\
&&+\int^s_{t_0}\int_E\beta(\tau^{-1}(r),
\wt{X}_{r-}^1,u^1_{\tau^{-1}(r)},e)\wt{\tau(\mu)}^{\mathbb{Q}_{\tau}}(\d
\tau^{-1}(r), \d e)\nonumber\\
&=&x_1+\int^s_{t_0}\frac{1}{\dot{\tau}}
b(\tau^{-1}(r),\wt{X}_r^1,u^0_r)\d
r+\int^s_{t_0}\frac{1}{\sqrt{\dot{\tau}}}\sigma(\tau^{-1}(r),\wt{X}_r^1,u^0_r)\d
B_r\nonumber\\
&&+\int^s_{t_0}\int_E\beta(\tau^{-1}(r),\wt{X}_{r-}^1,u^0_r,e)\lt(\tilde{\mu}(\d
r, \d e)+\lt(1-\frac{1}{\dot{\tau}}\rt)\Pi(\d e)\d r\rt).
\end{eqnarray}
This time change in equation (\ref{eq:sdej2}) makes the processes $X^0$
and $\wt{X}^1$ comparable. More precisely, we have
\begin{lemma}\label{lemma_X0-X1}
There exists some constant $C_\delta$, only depending on the bounds
of $\sigma, b, \beta$, their Lipschitz constants, as well as on $\Pi(E)$ and
$\delta$, such that, for all $t\in[t_0,T]$,
\[
\E\lt[\sup_{t\le s\le T}\lt|X_s^0-\wt{X}_s^1\rt|^2\big|\F_t\rt]\le
C_\delta\lt(|t_0-t_1|^2+|X^0_t-\wt{X}_t^1|^2\rt).
\]
\end{lemma}
For the proof of this lemma we need the following estimates gotten
by an elementary straightforward computation (see \cite{BCQ}, \cite{BHL}).
\begin{lemma}\label{lemma_tau1}
There is a constant $C_\delta$ only depending on $\delta>0$,
such that for all $r\in[t_0,T]$, we have
$$\lt|1-\frac{1}{\dot{\tau}}\rt|+\lt|\tau^{-1}(r)-r\rt|
+\lt|1-\sqrt{\dot{\tau}}\rt|\le C_\delta
|t_0-t_1|.$$
\end{lemma}

\noindent\textit{Proof (of Lemma \ref{lemma_X0-X1}):}
By taking the difference between the SDEs
(\ref{eq:sdej1}) and (\ref{eq:sdej3}) and
after the conditional expectation of the supremum of
its square, we get from Lemma \ref{lemma_tau1} and the assumptions on the coefficients, for $s\in[t,T]$,
\begin{eqnarray*}
&&\E\lt[\sup_{r\in[t,s]}\lt|X_r^0-\wt{X}_r^1\rt|^2\Big|\F_t\rt]\nonumber\\
&=&\E\bigg[\bigg(\lt|X^0_t-\wt{X}_t^1\rt|+\int^s_{t}\lt|b(v,X_v^0,u^0_v)
-\frac{1}{\dot{\tau}} b(\tau^{-1}(v),\wt{X}_v^1,u^0_v)\rt|\d v\nonumber\\
&&\quad\quad+\int^s_{t}\int_E\lt|\beta(\tau^{-1}(v),\wt{X}_{v-}^1,u^0_v,e)\rt|
\lt|1-\frac{1}{\dot{\tau}}\rt|\Pi(\d e)\d v\nonumber\\
&&\quad\quad+\sup_{r\in[t,s]}\lt|\int^r_{t}
\lt(\sigma(v,X_v^0,u^0_v)-\frac{1}{\sqrt{\dot{\tau}}}
\sigma(\tau^{-1}(v),\wt{X}_v^1,u^0_v)\rt)\d
B_v\rt|\nonumber\\
&&\quad\quad+\sup_{r\in[t,s]}\lt|\int^r_{t}
\int_E\lt(\beta(v,X_{v-}^0,u^0_v,e)-
\beta(\tau^{-1}(v),\wt{X}_{v-}^1,u^0_v,e)\rt)\tilde{\mu}(\d
v, \d e)\rt|\bigg)^2\bigg|\F_t\bigg]\nonumber\\
&\le& C\lt|X^0_t-\wt{X}_t^1\rt|^2+
C\E\lt[\lt(\int^s_{t}\lt(\lt|1-\frac{1}{\dot{\tau}}\rt|
+\lt|\tau^{-1}(v)-v\rt|+
\lt|X_v^0-\wt{X}_v^1\rt|\rt)\d v\rt)^2\bigg|\F_t\rt]\nonumber\\
&&\quad+C\E\lt[\int^s_{t}\lt(\lt|1-
\frac{1}{\sqrt{\dot{\tau}}}\rt|+\lt|\tau^{-1}(v)-v\rt|+
\lt|X_v^0-\wt{X}_v^1\rt|\rt)^2\d v\bigg|\F_t\rt]\nonumber\\
&\le &C_\delta\lt(\lt|X_t^0-\wt{X}_t^1\rt|^2+|t_0-t_1|^2\rt)
+C_\delta\int^s_{t}
\E\lt[\lt|X_v^0-\wt{X}_v^1\rt|^2|\F_t\rt]\d v.\nonumber
\end{eqnarray*}
Finally, from Gronwall's
inequality, we have
\[
\E\lt[\sup_{s\in[t,T]}\lt|X_s^0-\wt{X}_s^1\rt|^2\big|\F_t\rt]\le
C_\delta\lt(|t_0-t_1|^2+|X^0_t-\wt{X}_t^1|^2\rt).
\]
Hence the proof of Lemma \ref{lemma_X0-X1} is complete now.
\hfill$\cajita$
\bigskip

After having made comparable $X^0$ and $X^1$ by the time change of $X^1$,
we make now $Y^0$ and $Y^1$ comparable. For this we put
$\wt{Y}_s^1=Y_{\tau^{-1}(s)}^1$,
$\wt{Z}_s^1=\frac{1}{\sqrt{\dot{\tau}}}Z_{\tau^{-1}(s)}^1$
and $\wt{U}_s^1=U_{\tau^{-1}(s)}^1,$ $s\in[t_0,T]$.
Then $(\wt{Y}^1, \wt{Z}^1, \wt{U}^1)=(\wt{Y}^1_s,
\wt{Z}^1_s, \wt{U}^1_s)_{s\in[t_0,T]}\in
\mathcal{S}^2(t_0,T;\R)\times \mathbb{L}^2(t_0,T;\R^d)\times \mathbb{L}^2(t_0,T;\tilde{\mu},\R)$
is the solution of the BSDE
\begin{equation}\label{eq:bsdej3}
\begin{aligned}
\wt{Y}^1_s=&\Phi(\wt{X}_T^1)+\int^T_{s}\frac{1}{\dot{\tau}}
f(\tau^{-1}(r),\wt{X}_r^1,\wt{Y}_r^1,\sqrt{\dot{\tau}}\wt{Z}_r^1,
\wt{U}_r^1,u^0_r)\d r-\int^T_s
\wt{Z}_r^1\d B_r\\&-\int^T_s\int_E
\wt{U}^1_r\lt(\tilde{\mu}(\d r,\d
e)+\lt(1-\frac{1}{\dot{\tau}}\rt)\Pi(\d e)\d r\rt), \ s\in[t_0,T],
\end{aligned}
\end{equation}
with respect to the same filtration $\mathbb{F}$ as $(Y^0,Z^0,U^0)$.%

For the above BSDE, we have the following a priori estimates which can be proven by a straight-forward standard argument:
\begin{lemma}
Under hypothesis $(\mathbf{H2})$, there exists some constant $C_\delta$, only depending on the bounds
of $\sigma, b, \beta$, their Lipschitz constants, as well as on $\Pi(E)$ and
$\delta$, such that,
\begin{equation}\label{eq:apri}
\begin{aligned}
&\E\lt[\sup_{s\in[t,T]}|\wt{Y}^1_s|^2+
\int^T_t\lt|\wt{Z}_r^{1}\rt|^2\d r+
\int^T_t\int_E\lt|\wt{U}_r^{1}(e)\rt|^2\Pi(\d e)\d r\bigg|\F_t\rt]\le C_\delta<+\infty, \ \ t\in[t_0,T].
\end{aligned}
\end{equation}
\end{lemma}
Now we can state the key lemma for proving the joint Lipschitz continuity of $V$.
\begin{lemma}
Under our standard assumptions $(\mathbf{H1})$ and $(\mathbf{H2})$, we have
\begin{equation}\label{eq:est_BSDE}
\begin{aligned}
&\E\lt[\sup_{s\in[t,T]}|{Y}^0_s-\wt{Y}^1_s|^2+
\int^T_t\lt|{Z}_r^{0}-\wt{Z}_r^{1}\rt|^2\d r+
\int^T_t\int_E\lt|{U}_r^{0}(e)-\wt{U}_r^{1}(e)\rt|^2\Pi(\d e)\d r\bigg|\F_t\rt]\\
\le& C_\delta\lt(|t_0-t_1|^2+|X^0_t-\wt{X}_t^1|^2\rt), \ \ t\in[t_0,T].
\end{aligned}
\end{equation}
\end{lemma}
\noindent\textit{Proof}:
First we notice that, for $s\ge t$,
\begin{equation*}
\begin{aligned}
{Y}^{0}_s-\wt{Y}^1_s
=&\Phi({X}_T^{0})-\Phi(\wt{X}_T^{1})\\
&+\int^T_s\lt(
f(r,{X}_r^{0},
{Y}_r^{0},{Z}_r^{0},
{U}_r^{0},u^\lambda_r)-\frac{1}{\dot{\tau}}f(\tau^{-1}_1(r),\wt{X}_r^{1},
\wt{Y}_r^{1},\sqrt{\dot{\tau}}\wt{Z}_r^{1},
\wt{U}_r^{1},u^\lambda_r)\rt)\d
r\\&-\int^T_s \lt({Z}_r^{0}-\wt{Z}_r^{1}\rt)\d B_r-
\int^T_s\int_E\lt({U}_r^{0}(e)-\wt{U}_r^{1}(e)\rt)\widetilde{\mu}(\d
r, \d e)\\
&+ \int^T_s\int_E\lt(1-\frac{1}{\dot{\tau}}\rt)\wt{U}_r^{1}(e)\d r\Pi (\d e).
\end{aligned}
\end{equation*}
We apply It\^o's formula to $|{Y}^{0}_s-\wt{Y}^1_s|^2$ and,
using the boundedness and the Lipschitz continuity of $\Phi$ and $f$,
as well as Lemma \ref{lemma_tau1}, we deduce that
\begin{eqnarray*}
&&|{Y}^0_s-\wt{Y}^1_s|^2+\int^T_s\lt|{Z}_r^{0}-\wt{Z}_r^{1}\rt|^2\d r
+\int^T_s\int_E\lt|{U}_r^{0}(e)-\wt{U}_r^{1}(e)\rt|^2\Pi(\d e)\d r\nonumber\\
&\le&\lt|\Phi({X}_T^{0})-\Phi(\wt{X}_T^{1})\rt|^2
+C\int^T_s\lt|{X}_r^{0}-\wt{X}_r^{1}\rt|^2\d r
+C\int^T_s|{Y}^0_r-\wt{Y}^1_r|^2\d r\nonumber\\
&&-2\int^T_s\lt({Y}^0_r-\wt{Y}^1_r\rt)\lt({Z}^0_r-\wt{Z}^1_r\rt)\d B_r
+C|t_0-t_1|^2\nonumber\\
&&-\int^T_s\int_E\lt(2\lt({Y}^0_r-\wt{Y}^1_r\rt)\lt({U}^0_r(e)-\wt{U}^1_r(e)\rt)
+\lt|{U}_r^{0}(e)-\wt{U}_r^{1}(e)\rt|^2\rt)\tilde{\mu}(\d r,\d e)\nonumber\\
&&+C|t_0-t_1|^2\int^T_s\lt(\lt|\wt{Z}_r^1\rt|^2+
\int_E\lt|\wt{U}_r^{1}(e)\rt|^2\Pi(\d e)\rt)\d r.\nonumber
\end{eqnarray*}
By taking the conditional expectation on both sides,
using Lemma \ref{lemma_X0-X1}, the a priori estimate (\ref{eq:apri}) and Gronwall's lemma,
we obtain, for $t_0\le t\le s\le T$,
\[
\begin{aligned}
&\E\lt[|{Y}^0_s-\wt{Y}^1_s|^2+\int^T_s\lt|{Z}_r^{0}-\wt{Z}_r^{1}\rt|^2\d r
+\int^T_s\int_E\lt|{U}_r^{0}(e)-\wt{U}_r^{1}(e)\rt|^2\Pi(\d e)\d r\bigg|\F_t\rt]\\
\le& C_\delta\lt(|t_0-t_1|^2+|X^0_t-\wt{X}_t^1|^2\rt).
\end{aligned}
\]
Then the Burkholder-Davis-Gundy inequality allows to show that
\[
\begin{aligned}
&\E\lt[\sup_{s\in[t,T]}|{Y}^0_s-\wt{Y}^1_s|^2+
\int^T_t\lt|{Z}_r^{0}-\wt{Z}_r^{1}\rt|^2\d r+
\int^T_t\int_E\lt|{U}_r^{0}(e)-\wt{U}_r^{1}(e)\rt|^2\Pi(\d e)\d r\bigg|\F_t\rt]\\
\le& C_\delta\lt(|t_0-t_1|^2+|X^0_t-\wt{X}_t^1|^2\rt),\ \ t\in[t_0,T].
\end{aligned}
\]
The proof is now complete.
\hfill$\cajita$

Now we are ready to give the proof of Theorem \ref{thm:lip}.
\bigskip

\noindent\textit{Proof of Theorem \ref{thm:lip}}:
By taking $t=t_0$ in (\ref{eq:est_BSDE}), we have
$$
\begin{aligned}
&|Y^0_{t_0}-Y^1_{t_1}|^2=|{Y}^0_{t_0}-\wt{Y}^1_{t_0}|^2\\
\le& C_\delta\lt(|t_0-t_1|^2+|X^0_{t_0}-\wt{X}_{t_0}^1|^2\rt)
=C_\delta\lt(|t_0-t_1|^2+|X^0_{t_0}-X^1_{t_1}|^2\rt)\\
=&C_\delta\lt(|t_0-t_1|^2+|x_0-x_1|^2\rt).\\
\end{aligned}
$$
Therefore, from (\ref{cost1}) and (\ref{cost2}), we get that
\[
\begin{aligned}
&V(t_0,x_0)-V(t_1,x_1)\\
\ge& J(t_0,x_0;u^0)-J(t_1,x_1;u^1)-\ep\\
=&Y^0_{t_0}-\wt{Y}^1_{t_0}-\ep\\
\ge& -C_\delta(|t_0-t_1|+|x_0-x_1|)-\ep,
\end{aligned}
\]
for some $C_\delta$ only depending on $\delta$ but not on $(t_0,x_0),$
$(t_1,x_1)\in[0,T-\delta]\times\R^d$.
Thus, from the arbitrariness of $\ep$, we deduce that 
\[
V(t_0,x_0)-V(t_1,x_1)\ge -C_\delta(|t_0-t_1|+|x_0-x_1|).
\] 
Symmetrical argument yields the converse relation. Consequently, the joint Lipschitz
continuity of $V$ over $[0,T-\delta]\times\R^d$.\hfill$\cajita$

\section{Semiconcavity}
\setcounter{equation}{0}
We study in this section the semiconcavity property of
the viscosity solution $V$ and to extend for this the method of
time change and Kulik's transformation used in the preceding section.

For the semiconcavity property,
we need more assumptions on the coefficients:

\noindent\textbf{(H3)} The function $\Phi(x)$ is
semiconcave, and $f(\cdot,\cdot,\cdot,\cdot,\cdot,u)$ is semiconcave in
$(t,x,y,z,p)\in[0,T]\times\R^d\times\R\times\R^d\times L^2(E,\mathcal{B}(E),\Pi;\R)$,
uniformly with respect to $u\in U$, i.e.,
there exists a constant $C>0$, such that, for any
$\xi_1\triangleq(t_1,x_1,y_1,z_1,p_1)$,
$\xi_2\triangleq(t_2,x_2,y_2,z_2,p_2)$
in $[0,T]\times\R^d\times\R\times\R^d\times L^2(E,\mathcal{B}(E),\Pi;\R)$, and $\lambda\in[0,1]$,  $u\in U$,
\[
\begin{aligned}
&\lambda f(\xi_1,u)+(1-\lambda)f(\xi_2,u)-f(\lambda \xi_1+(1-\lambda)\xi_2,u)\\
\le& C\lambda(1-\lambda)\lt(|t_1-t_2|^2+|x_1-x_2|^2+|y_1-y_2|^2+|z_1-z_2|^2+\int_E|p_1(e)-p_2(e)|^2\Pi(\d e)\rt).
\end{aligned}
\]

\noindent\textbf{(H4)} The first-order derivatives $\nabla_{t,x}b$,
$\nabla_{t,x}\sigma$ and $\nabla_{t,x}\beta$ of $b,$ $\sigma$ and
$\beta$ with respect to $(t,x)$ exist and are continuous in $(t,x,u)$ and Lipschitz
continuous in $(t,x)$, uniformly with respect to $u$.

\noindent\textbf{(H5)}  There exist two constants $-1<C_1<0$ and  $C_2>0$
such that, for all $(t,\xi):=(t,x,y,z)\in[0,T]\times\R^d\times\R\times\R^d$,
$u\in U$, and $p, p^\pr\in L^2(E,\mathcal{B}(E),\Pi;\R)$,
$$f(t,\xi,p,u)-f(t,\xi,p^\pr,u)\le \int_E(p(e)-p^\pr(e))\gamma_t^{\xi,u;p,p^\pr}(e)\Pi(\d e),$$
where $\{\gamma_t^{\xi,u;p,p^\pr}(e)\}_{t\in[0,T]}$ is a measurable function such that,
for every $t\in[0,T]$,
$$C_1(1\wedge |e|)\le \gamma_t^{\xi,u;p,p^\pr}(e)\le C_2(1\wedge |e|).$$

Our main result in this section is the following:
\begin{theorem}\label{thm:semi}
Under the assumptions $(\textbf{H1})-(\textbf{H5})$, for every $\delta\in (0,T)$, there
exists some constant $C_\delta>0$ such that, for all $(t_0,x_0),
(t_1,x_1)\in[0,T-\delta]\times\R^d$, and for all $\lambda\in[0,1]$:
\[
\lambda V(t_0,x_0)+(1-\lambda)V(t_1,x_1)-V(t_\lambda,x_\lambda)\le
C_\delta\lambda(1-\lambda)(|t_0-t_1|^2+|x_0-x_1|^2),
\]
where $t_\lambda=\lambda t_0+(1-\lambda)t_1, x_\lambda=\lambda
x_0+(1-\lambda)x_1.$
\end{theorem}

\begin{remark}
Again as in the case of the Lipschitz continuity,
we cannot hope, in general, that the property of semiconcavity holds
over the whole domain $[0,T]\times\R^d$.
Indeed, let us consider the example given in the preceding section (Remark \ref{remark_lip}). 
In particular, we have gotten there that 
$$
V(s,0)=\E[\Phi(X_T^{s,0})]=\E[-|B_T-B_s|]=-\sqrt{\frac{2}{\pi}}\sqrt{T-s},\ \ \ s\in[0,T].
$$
However, it is easy to check that this function $V$ is not semiconcave in $[0,T]\times\R^d$, 
but it has this semiconcavity property on $[0,T-\delta]\times\R^d$, for all $\delta>0$.
\end{remark}

The proof of  Theorem \ref{thm:semi} will be based again on the method of time change.
But unlike the proof of the Lipschitz property, we have to work here with two time changes.
In order to be more precise, for given
$\delta>0$, $(t_0,x_0), (t_1,x_1)\in[0,T-\delta]\times\R^d$,
let us consider the both following linear time changes:
\[
\tau_i: [t_i,T]\to [t_\lambda,T],\
\tau_i(t)=t_\lambda+\frac{T-t_\lambda}{T-t_i}(t-t_i),\]
with the derivatives
$\dot{\tau}_i=\frac{T-t_\lambda}{T-t_i}, \ t\in [t_i,T], \ i=0,1.$

For $t=t_{\lambda}$, we let $B=\{B_s\}_{s\in[t_\lambda,T]}$ be a Brownian motion
starting from zero at $t_\lambda$: $B_{t_{\lambda}}=0$. Then
$\{W^i_s=\frac{1}{\sqrt{\dot{\tau_i}}}B_{\tau_i(s)}\}_{s\in[t_i,T]}$ is a
Brownian motion on $[t_i,T]$, starting from zero at time $t_i$,
$i=0,1$. For $t=t_\lambda$,  let $\mu(\d r, \d e)$ be our Poisson random measure on
$[t_\lambda,T]\times E$ under probability $\PP$. Then $\tau_i(\mu), i=0,1$,
defined as the Kulik transformation of $\mu$,
\[
\tau_i(\mu)([t_i,t_i+s]\times\Delta)\triangleq\mu([t_\lambda,
\tau_i(t_i+s)]\times\Delta), 0\le s\le T-t_i, \Delta\in \mathcal{B}E
\]
is a new Poisson random measure but under probability $\mathbb{Q}_i$, where
\[
\frac{\d \mathbb{Q}_{i}}{\d \PP}=\exp\lt\{-\ln
\lt(\frac{T-t_\lambda}{T-t_i}\rt)\mu([t_i,T]\times
E)+(t_i-t_\lambda)\Pi(E)\rt\}.
\]
We denote the corresponding compensated Poisson random measures
under $\PP$ and $\mathbb{Q}_i$ by $\tilde{\mu}$
and $\wt{\tau_i(\mu)},\ i=0,1,$ respectively:
\[
\tilde{\mu}(\d s, \d e)={\mu}(\d s, \d e)-\d s \Pi(\d e),\
(s,e)\in[t_\lambda,T]\times E,
\]
and
\[
\wt{\tau_i(\mu)}(\d s, \d e)=\tau_i(\mu)(\d s, \d e)-\d s \Pi(\d e),\
(s,e)\in[t_i,T]\times E.
\]

Let us now fix an arbitrary  $u^\lambda\in \mathcal{U}^{B,\mu}(t_\lambda,T)$
(Recall the definition of $\mathcal{U}^{B,\mu}(t_\lambda,T)$).
Then, obviously, $u^i_s\triangleq u^\lambda_{\tau_i(s)}, s\in[t_i,T]$,
is an admissible control in
$\mathcal{U}^{W^i,\tau_i(\mu)}(t_i,T)$ with respect to $W^i$ and $\tau_i(\mu)$.

We let $\{X_s^{\lambda}\}_{s\in[t_\lambda, T]}$ be the unique solution of
the SDE,
\begin{equation}\label{eq:2.1}
X_s^{\lambda}=x_\lambda+\int^s_{t_\lambda}
{\hskip-1mm}b(r,X_r^{\lambda},u^\lambda_r)\d r +
\int^s_{t_\lambda}{\hskip-1mm}\sigma(r,X_r^{\lambda},u^\lambda_r)\d
B_r +\int^s_{t_\lambda}{\hskip-1mm}\int_E{\hskip-1mm}
\beta(r,X_{r-}^{\lambda},u^\lambda_r,e)\tilde{\mu}(\d
r, \d e),\ s\in[t_\lambda,T].
\end{equation}
We also make use of the unique solution
$\{X_s^i\}_{s\in[t_i,T]}$ of the following SDE,
\begin{equation}\label{eq:2.2}
X_s^i=x_i+\int^s_{t_i} {\hskip-1mm}b(r,X_r^i,u^i_r)\d r
+ \int^s_{t_i}{\hskip-1mm}\sigma(r,X_r^i,u^i_r)\d W_r^i
+\int^s_{t_i}{\hskip-1mm}\int_E{\hskip-1mm}
\beta(r,X_{r-}^i,u^i_{r},e)\wt{\tau_i(\mu)}(\d r, \d
e),\ s\in[t_i,T],  \ i=0,1.
\end{equation}

As in the preceding section, we associate the forward equations
(\ref{eq:2.1}) and (\ref{eq:2.2}) with BSDEs.
Let $(Y_s^{\lambda},
Z_s^{\lambda},
U_s^{\lambda})_{s\in[t_\lambda,T]}$ and
$(Y_s^i,Z_s^i,U_s^i)_{s\in[t_i,T]}$,
$i=0,1$, be the unique solutions of the BSDEs
\begin{eqnarray*}
Y_s^{\lambda}=\Phi(X_T^{\lambda})
+\int^T_s{\hskip-1mm}
f(r,X_r^{\lambda},Y_r^{\lambda},
Z_r^{\lambda},
U_r^{\lambda},u^\lambda_r)\d r -\int^T_s
{\hskip-1mm} Z_r^{\lambda}\d B_r-\int^T_s{\hskip-2mm}
\int_E{\hskip-1mm}U_r^{\lambda}(e)\tilde{\mu}(\d
r, \d e),\ s\in[t_\lambda,T],
\end{eqnarray*}
and
\begin{eqnarray*}
Y_s^{i}=\Phi(X_T^{i}) +\int^T_s{\hskip-1mm}
f(r,X_r^{i},Y_r^i,Z_r^i,U_r^i,u^i_r)\d
r -\int^T_s {\hskip-1mm}Z_r^{i}\d W_r^i-\int^T_s{\hskip-2mm}
\int_E{\hskip-1mm}U_r^{i}(e)\widetilde{\tau_i(\mu)}(\d r,
\d e),\ s\in[t_i,T],
\end{eqnarray*}
respectively. Then from the adaptedness of the solutions $Y^\lambda$
and $Y^i$ with respect to the filtrations generated by $(B,\mu)$
and $(W^i,\tau_i(\mu))$, respectively, we know
that $Y_{t_\lambda}^{\lambda}$ and
$Y_{t_i}^{i}$ are deterministic and equal to the cost functionals
$J(t_\lambda,x_\lambda;u^\lambda)$ and $J(t_i,x_i;u^i)$,
respectively.

For the proof of the Theorem \ref{thm:semi} it is crucial to estimate
$\lambda Y^0_{t_0}+(1-\lambda)Y_{t_1}^1-Y_{t_\lambda}^\lambda,$ for $\lambda\in(0,1).$
Since the processes $Y^0, Y^1$ and $Y^\lambda$ are solutions of
BSDEs over different time intervals, driven by different Brownian motions
and different Poisson random measures, we have to make them comparable
with the help of the inverse time change.

In a first step we carry out this inverse time change for the forward equations.
For this end we introduce the time-changed processes:
$\widetilde{X}_s^i=X_{\tau_i^{-1}(s)}^i$, $s\in[t_\lambda,T],$ $i=0,1$.
Then we have, for $i=0,1,$
\begin{equation}
\begin{aligned}
\widetilde{X}_s^i=&x_i+\int^s_{t_\lambda}
\frac{1}{\dot{\tau_i}}b(\tau^{-1}_i(r),\widetilde{X}_r^i,u^\lambda_r)\d
r + \int^s_{t_\lambda}\frac{1}{\sqrt{\dot{\tau_i}}}\
\sigma(\tau^{-1}_i(r),\widetilde{X}_r^i,u^\lambda_r)\d
B_r \\
&+\int^s_{t_\lambda}\int_E
\beta(\tau^{-1}_i(r),\widetilde{X}_{r-}^i,u^\lambda_r,e)
\lt(\widetilde{\mu}(\d r, \d e)+\lt(1-\frac{1}{\dot{\tau_i}}\rt)\Pi (\d e)\d
r\rt), \ s\in[t_\lambda,T].
\end{aligned}
\end{equation}

Comparable with Lemma \ref{lemma_X0-X1} but now with arbitrary power $p\ge 2$,
we can show
\begin{lemma}\label{lemma_1-2} Let $p\ge 2$.
Then there exists a constant $C_{\delta,p}$ depending
only on the bounds of $\sigma, b, \beta$, their Lipschitz constants,
$\Pi(E)$, $\delta$ as well as $p$, such that, for all $t\in[t_\lambda,T]$,
\begin{equation}\label{eq:^2}
\E\lt[\sup_{t\le s\le T}\lt|\wt{X}_s^{0}
-\wt{X}_s^{1}\rt|^p\bigg|\F_t\rt]\le C_{\delta,p}
\lt(\lt|\wt{X}_t^{0}-\wt{X}_t^{1}\rt|^p+
\lt|t_0-t_1\rt|^p\rt).
\end{equation}
\end{lemma}
Moreover, in addition to Lemma \ref{lemma_1-2}, which gives a kind of
``first order estimate'', we also have the following kind of ``second order estimate''.
For this we introduce the process
$\wt{X}^{\lambda}=\lambda\wt{X}^{0}+(1-\lambda)\wt{X}^{1}$.
\begin{lemma}\label{lemma_lam-lam}
Let $p\ge 2$. There exists a constant $C_{p,\delta}$ depending
only on the bounds of $\sigma, b, \beta$, their Lipschitz constants,
$\Pi(E)$, $\delta$ and $p$, such that, for all $t\in[t_\lambda,T]$,
\begin{equation}\label{eq:lam}
\begin{aligned}
&\E\lt[\sup_{t\le s\le T}
\lt|\wt{X}_s^{\lambda}-{X}_s^{{\lambda}}\rt|^p\bigg|\F_t\rt]\\
\le&
C_{p,\delta}\lt|\wt{X}_t^{\lambda}-{X}_t^{{\lambda}}
\rt|^p+C_{p,\delta}(\lambda(1-\lambda))^{p}\lt(|t_0-t_1|^{2p}+\lt|\wt{X}_t^{0}
-\wt{X}_t^{1}\rt|^{2p}\rt).
\end{aligned}
\end{equation}
\end{lemma}
For the proof of the both above lemmata the reader is referred to the Appendix.
\bigskip

After having applied the inverse time changes to the forward equations,
let us do it now for the BSDEs. Thus, for $i=0,1$, we introduce the processes
$\wt{Y}^i_s\triangleq
Y^i_{\tau^{-1}_i(s)}$, $\wt{Z}_s^{i}\triangleq\frac{1}{\sqrt{\dot{\tau_i}}}
Z_{\tau^{-1}_i(s)}^{i}$,
and $\wt{U}_s^{i}\triangleq U_{\tau^{-1}_i(s)}^{i}$,  $s\in[t_\lambda,T]$.
Obviously, $\lt(Y^\lambda,Z^\lambda,U^\lambda\rt)$
and $(\wt{Y}^i,\wt{Z}^i,\wt{U}^i)$ belong to
$\mathcal{S}^2(t_\lambda,T;\R)\times \mathbb{L}^2(t_\lambda,T;\R^d)\times
\mathbb{L}^2(t_\lambda,T;\tilde{\mu},\R)$, and
\begin{equation*}
\begin{aligned}
\wt{Y}^i_s=&\Phi(\wt{X}_T^{i}) +\int^T_s
\frac{1}{\dot{\tau}_i}f(\tau^{-1}_i(r),\wt{X}_r^{i},
\wt{Y}_r^i,\sqrt{\dot{\tau}_i}\wt{Z}_r^i,
\wt{U}_r^i,u^\lambda_r)\d
r -\int^T_s \wt{Z}_r^{i}\d B_r\\
&\quad-\int^T_s\int_E\wt{U}_r^{i}(e)\lt(\widetilde{\mu}(\d
r, \d e)+\lt(1-\frac{1}{\dot{\tau_i}}\rt)\Pi (\d e)\d r\rt),
 \ s\in[t_\lambda,T].
\end{aligned}
\end{equation*}
With the help of standard BSDE estimates we can show
\begin{lemma} \label{lemma_apriori} For $p\ge 2$,
there exists some constant $C_p$ only depending on $p$
and the bounds of the coefficients $f, \Phi$,
such that, for all $s\in[t_\lambda,T]$, $i=0, 1$,
\begin{equation*}
\E\lt[\sup_{s\le r\le T}|\wt{Y}^i|^p+
\lt(\int^T_s|\wt{Z}^i_r|^2\d r\rt)^{p/2}
+\lt(\int^T_s\int_E|\wt{U}^i_r(e)|^2\Pi(\d e)\d r\rt)^{p/2}\bigg|\F_s\rt]\le C_p,
\end{equation*}
and
\begin{equation*}
\E\lt[\sup_{s\le r\le T}|{Y}^\lambda|^p
+\lt(\int^T_s|{Z}^\lambda_r|^2\d r\rt)^{p/2}
+\lt(\int^T_s\int_E|{U}^\lambda_r(e)|^2\Pi(\d e)\d r\rt)^{p/2}\bigg|\F_s\rt]\le C_p.
\end{equation*}
\end{lemma}
As the proof uses simple BSDE estimates
which by now are standard (see, for instance, \cite{CT}), the proof is omitted.

Recall that we have defined
$\wt{X}^{\lambda}=\lambda\wt{X}^{0}+(1-\lambda)\wt{X}^{1}$.
In the same manner, we introduce the processes
$\wt{Y}^{\lambda}=\lambda\wt{Y}^{0}+(1-\lambda)\wt{Y}^{1}$,
$\wt{Z}^{\lambda}=\lambda\wt{Z}^{0}+(1-\lambda)\wt{Z}^{1}$ and
$\wt{U}^{\lambda}=\lambda\wt{U}^{0}+(1-\lambda)\wt{U}^{1}$. Then we
get that
$(\wt{Y}^\lambda, \wt{Z}^\lambda, \wt{U}^\lambda)
\in \mathcal{S}^2(t_\lambda,T;\R)\times
\mathbb{L}^2(t_\lambda,T;\R^d)\times
\mathbb{L}^2(t_\lambda,T;\tilde{\mu},\R)$
is the unique solution of the following BSDE
\[
\begin{aligned}
&\wt{Y}_s^\lambda=\lambda\Phi\lt(\wt{X}_T^{0}\rt)+
(1-\lambda)\Phi\lt(\wt{X}_T^{1}\rt)-\int^T_s \wt{Z}_r^\lambda\d
B_r-\int^T_s\int_E\wt{U}_r^{\lambda}(e)\tilde{\mu}(\d
r, \d
e)\\
&+\int^T_s\hspace{-1mm}
\lt[\frac{\lambda}{\dot{\tau_0}}f(\tau_0^{-1}(r),
\wt{X}_r^0,\wt{Y}_r^0,\sqrt{\dot{\tau_0}}\wt{Z}_r^0,\wt{U}_r^0,u^\lambda_r)
+\frac{1-\lambda}{\dot{\tau_1}}f(\tau_1^{-1}(r),
\wt{X}_r^1,\wt{Y}_r^1,\sqrt{\dot{\tau_1}}\wt{Z}_r^1,\wt{U}_r^1,u^\lambda_r)\rt]\d
r\\
&-\int^T_s\int_E\lt[\lambda\lt(1-\frac{1}{\dot{\tau_0}}\rt)\wt{U}_r^{0}(e)+
(1-\lambda)\lt(1-\frac{1}{\dot{\tau_1}}\rt)\wt{U}_r^{1}(e)\rt]\Pi (\d e)\d
r, \  s\in[t_\lambda,T].
\end{aligned}
\]

In analogy to Lemma \ref{lemma_1-2}
we have for the associated BSDE the following statement,
which proof is postponed in the Appendix:
\begin{lemma}\label{lemma_Y0-Y1}
For all $p\ge 2$, there exists a constant $C_\delta$ depending
only on the bounds of $\sigma, b, \beta$, their Lipschitz constants,
$\Pi(E)$, $\delta$ and $p$, such that, for any $t\in[t_\lambda,T]$,
\begin{equation*}
\begin{aligned}
&\E\lt[\sup_{t\le s\le T}|\wt{Y}^0_s-\wt{Y}^1_s|^p
+\lt(\int^T_t\lt|\wt{Z}^0_s-\wt{Z}^1_s\rt|^2\d s\rt)^{p/2}
+ \lt(\int^T_t\int_E\lt|\wt{U}^0_s(e)-\wt{U}^1_s(e)\rt|^2\Pi(\d e)
\d s\rt)^{p/2}\bigg|\F_t\rt]\\
&\le C_\delta\lt(|\wt{X}_t^{0}
-\wt{X}_t^{1}|^p+|t_0-t_1|^p\rt).
\end{aligned}
\end{equation*}
\end{lemma}

Our objective is to estimate
\[
\lambda Y^0_{t_0}+(1-\lambda) Y^1_{t_1}-Y^{\lambda}_{t_{\lambda}}
=\wt{Y}_{t_{\lambda}}^{\lambda}-Y^{\lambda}_{t_\lambda}.
\]
For this end some auxiliary processes shall be introduced.
So let us introduce the increasing c\`adl\`ag processes
$$
A_t:= |t_0-t_1|+\sup_{s\in[t_\lambda,t]} |\wt{X}_s^0-\wt{X}_s^1|,
$$
and
$$
B_t:=\sup_{s\in[t_\lambda,t]}|\wt{X}_s^\lambda-X_s^\lambda|,\ \
t\in[t_\lambda,T].
$$
For some suitable $C$ and $C_\delta$
which will be specified later, we also introduce the increasing c\`adl\`ag process
$$
D_t=CB_t+C_\delta\lambda(1-\lambda)A_t^2,\ \ t\in[t_\lambda,T].
$$

We can obtain easily from the Lemmata \ref{lemma_1-2} and \ref{lemma_lam-lam}
the following estimate for $D_t$.
\begin{corollary}\label{cor_D}
For any $p\ge2$,
there exists a constant $C_p$ such that
\[
\E\lt[|D_s|^p|\F_t\rt]\le C_p |D_t|^p, \ \ \mbox{for\ all}
\ \ s,t\in[t_\lambda,T],\ \ \mbox{with}\ \ t\le s.
\]
\end{corollary}
We observe that, in particular,
$$
\E\lt[|D_T|^p\rt]\le C_p(|t_0-t_1|^p+|x_0-x_1|^p)<+\infty, \ \ p\ge2.
$$

We let
$(\wh{Y}^\lambda, \wh{Z}^\lambda, \wh{U}^\lambda)
\in \mathcal{S}^2(t_\lambda,T;\R)\times
\mathbb{L}^2(t_\lambda,T;\R^d)\times
\mathbb{L}^2(t_\lambda,T;\tilde{\mu},\R)$
be the unique solution of the following BSDE,
\begin{equation*}
\begin{aligned}
\wh{Y}_s^\lambda&=\Phi\lt({X}_T^{\lambda}\rt)+D_T-\int^T_s \wh{Z}_r^\lambda\d
B_r-\int^T_s\int_E\wh{U}_r^{\lambda}(e)\tilde{\mu}(\d
r, \d
e)\\
&+\int^T_s
\bigg[f(r,{X}_r^{\lambda},\wh{Y}_r^\lambda-D_r,\wh{Z}_r^\lambda,\wh{U}_r^\lambda,u^\lambda_r)
+CD_r\\&
+C_\delta^0\lambda(1-\lambda)\lt(|t_0-t_1|^2\lt(1+|\wt{Z}_r^1|^2\rt)+|\wt{Z}_r^0-\wt{Z}_r^1|^2\rt)
+\int_E\lt|\wt{U}_r^{0}(e)-\wt{U}_r^{1}(e)\rt|^2\Pi(\d e)\bigg]
\d r,\ \ s\in[t_\lambda,T].
\end{aligned}
\end{equation*}
The process $\wh{Y}^{\lambda}$ stems its importance
from the fact that it majorizes $\wt{Y}^{\lambda}$  in a suitable manner.
More precisely, we have
\begin{lemma}\label{lemma_comparison}
$\wt{Y}_s^\lambda\le\wh{Y}_s^\lambda$, $\PP$-a.s., for any $s\in[t_\lambda,T]$.
\end{lemma}
For a better readability of the paper, also this proof is postponed to the Appendix.
\bigskip

In addition to Lemma \ref{lemma_comparison}, we also have to
estimate the difference between $\wt{Y}^{\lambda}$ and $Y^\lambda$. For this we
introduce the process $\overline{Y}_t^\lambda=\wh{Y}_t^\lambda-D_t$, $t\in[t_\lambda,T]$,
and we identify $(\overline{Y}^\lambda,\wh{Z}^\lambda,\wh{U}^\lambda)$
as the unique solution of the BSDE
\begin{equation*}
\begin{aligned}
\overline{Y}_s^\lambda=&\Phi\lt({X}_T^{\lambda}\rt)+\int^T_s
\bigg[f(r,{X}_r^{\lambda},\overline{Y}_t^\lambda,\wh{Z}_r^\lambda,\wh{U}_r^\lambda,u^\lambda_r)
+CD_r\\
&\quad+C_\delta^0\lambda(1-\lambda)\lt(|t_0-t_1|^2\lt(1+|\wt{Z}_r^1|^2\rt)+|\wt{Z}_r^0-\wt{Z}_r^1|^2\rt)\\
&\quad
+\int_E\hspace{-1mm}\lt|\wt{U}_r^{0}(e)-\wt{U}_r^{1}(e)\rt|^2\Pi(\d e)\bigg]
\d r-\int^T_s \hspace{-1mm}\wh{Z}_r^\lambda\d
B_r-\int^T_s\hspace{-2mm}\int_E\wh{U}_r^{\lambda}(e)\tilde{\mu}(\d
r, \d e)+\int_s^T\hspace{-1mm}\d D_r,\ \ s\in[t_\lambda,T].
\end{aligned}
\end{equation*}
We observe that we have the following statement,
which proof is given in the Appendix.

\begin{lemma} \label{lemma_wtY-Y}
For $t\in[t_\lambda,T]$, we have
\begin{equation}\label{eq:_Y-Y}
\begin{aligned}
&\E\lt[\sup_{s\in[t,T]}\lt|\overline{Y}_s^\lambda-Y^\lambda_s\rt|^2
+\int^T_t\lt|\wh{Z}_s^\lambda-Z^\lambda_s\rt|^2\d s
+\int^T_t\int_E \lt|\wh{U}_s^{\lambda}(e)-U_s^\lambda(e)\rt|^2{\mu}(\d
s, \d e)\bigg|\F_t\rt]\\
\le &C_\delta D_t^2.
\end{aligned}
\end{equation}
\end{lemma}

Now we are ready to give the proof of Theorem \ref{thm:semi}.

\noindent\textit{Proof of Theorem \ref{thm:semi}}:
We know from the stochastic interpretation of the viscosity solution $V$ as value function (see (\ref{eq:value}))
%\[
%V(t_{\lambda},x_\lambda)=\inf_{u\in\mathcal{U}^{B,\mu}(t_\lambda,T)} Y_{t_\lambda}^\lambda
%\]
that, for any $\lambda\in(0,1)$ and $\ep>0$, there exists an admissible control process
$u^\lambda\in\mathcal{U}^{B,\mu}(t_\lambda,T)$ such that
$Y_{t_\lambda}^{\lambda}\le V(t_\lambda,x_\lambda)+\ep$.
On the other hand, using again (\ref{eq:value}), but now for $\mathcal{U}^{W^i,\tau_i(\mu)}$, we obtain:
$V(t_i,x_i)\le Y_{t_i}^i,$ $i=0,1.$
From the Lemmata \ref{lemma_comparison} and \ref{lemma_wtY-Y} we deduce that
\begin{eqnarray*}
&&\lambda V(t_0,x_0)+(1-\lambda)V(t_1,x_1)\\
&\le& \lambda {Y}^{0}_{t_0}+(1-\lambda) {Y}^{1}_{t_1}\\
&=&
\lambda\wt{Y}^{0}_{t_\lambda}+(1-\lambda)\wt{Y}^{1}_{t_\lambda}
=\wt{Y}^\lambda_{t_\lambda}\\
&\le& \wh{Y}_{t_\lambda}^\lambda=\overline{Y}_{t_\lambda}^{\lambda}+D_{t_\lambda}\\
&\le& Y^\lambda_{t_\lambda}+C_\delta D_{t_\lambda}\\
&\le&V(t_\lambda,x_\lambda)+C_\delta B_{t_\lambda}
+C_\delta \lambda(1-\lambda) A_{t_\lambda}^2+\ep\\
&\le& V(t_\lambda,x_\lambda)+ C_\delta \lambda(1-\lambda)
(|t_0-t_1|^2+|x_0-x_1|^2)+\ep.
\end{eqnarray*}
Here we have used that
$D_{t_\lambda}=CB_{t_\lambda}+C_\delta\lambda(1-\lambda)A^2_{t_\lambda}$
and $B_{t_\lambda}=0$.
From the arbitrariness  of $\ep$, it follows that
\[
\lambda V(t_0,x_0)+(1-\lambda)V(t_1,x_1)
\le V(t_\lambda,x_\lambda)+ C_\delta \lambda(1-\lambda) (|t_0-t_1|^2+|x_0-x_1|^2).
\]
Hence, the semiconcavity of $V$ is proved.\hfill$\cajita$

\section{Appendix}
\setcounter{equation}{0}
The appendix is devoted to the proof of the Lemmata \ref{lemma_1-2}--\ref{lemma_wtY-Y}.

First we give the following lemma, which will be used in what follows.
It can be checked by a straightforward computation and, hence, its proof is omitted.
\begin{lemma}\label{lemma_tau}\label{}
There exists some positive constant $C_\delta$ only depending on $T$
and $\delta$ such that, for $s\in[t_\lambda,T]$,
$$|\tau^{-1}_0(s)-\tau^{-1}_1(s)|+|\frac{1}{\dot{\tau}_0}-\frac{1}{\dot{\tau}_1}|
+|\frac{1}{\sqrt{\dot{\tau}_0}}-\frac{1}{\sqrt{\dot{\tau}_1}}|\le
C_\delta |t_0-t_1|,$$
\[
\lambda\lt|1-\frac{1}{\sqrt{\dot{\tau_0}}}\rt|+(1-\lambda)\lt|1-\frac{1}{\sqrt{\dot{\tau_1}}}\rt|
\le \frac{1}{2\delta}\lambda(1-\lambda)|t_0-t_1|,
\]
\[
\lt|\lambda(1-\frac{1}{\sqrt{\dot{\tau_0}}})+(1-\lambda)(1-\frac{1}{\sqrt{\dot{\tau_1}}})\rt|\le \frac{1}{\delta^2}\lambda(1-\lambda)|t_0-t_1|^2.
\]
Moreover, for all $s\in[t_\lambda,T],$
$$\lambda(1-\frac{1}{\dot{\tau}_0})=-(1-\lambda)(1-\frac{1}{\dot{\tau}_1})
=\frac{\lambda(1-\lambda)}{T-t_\lambda}(t_0-t_1),\
\lambda\tau^{-1}_0(s)+(1-\lambda)\tau^{-1}_1(s)=s.$$
\end{lemma}

We begin with the proof of Lemma \ref{lemma_1-2}.
\bigskip

\noindent\textit{Proof of Lemma \ref{lemma_1-2}:} Let us put
$$
\Delta \beta(r,e)=\beta(\tau^{-1}_0(r),
\widetilde{X}_{r-}^{0},u^\lambda_{r},e)-
\beta(\tau^{-1}_1(r),\widetilde{X}_{r-}^{1},
u^\lambda_{r},e),
$$
for $(r,e)\in[t_\lambda,T]\times E$,  and consider, for $t\in[t_\lambda,T]$,
$$
M_s=\int^s_{t}\int_E\Delta \beta(r,e)\tilde{\mu}(\d
r,\d e), \  s\in[t,T].
$$
Since $\Delta \beta$ is bounded and predictable,
$M$ is a $p$-integrable martingale, for all $p\ge 2$. We deduce from It\^o's formula that,
$$
N_s: =|M_s|^p-\int^s_t\int_E\lt(|M_r+\Delta\beta(r,e)|^p-|M_r|^p-p|M_r|^{p-2}M_r\Delta \beta(r,e)\rt)\Pi(\d e)\d r, \ s\in[t,T],
$$
is a martingale with $N_t=0$ (see, Fujiwara and Kunita \cite{FK}).
Moreover, since $\beta$ is Lipschitz in $(t,x)$, uniformly with respect to $(u,e)$,
$$
\begin{aligned}
&|M_r+\Delta \beta(r,e)|^p-|M_r|^p-p|M_r|^{p-2}M_r\Delta \beta(r,e)\\
\le& C_p\lt(|\Delta \beta(r,e)|^2|M_r|^{p-2}+|\Delta \beta(r,e)|^p\rt)\\
\le& C_p\lt(\lt(|\wt{X}_{r}^{0}
-\wt{X}_{r}^{1}|^2+|\tau_0^{-1}(r)-\tau_1^{-1}(r)|^2\rt)|M_r|^{p-2}+|\wt{X}_{r}^{0}
-\wt{X}_{r}^{1}|^p+|\tau_0^{-1}(r)-\tau_1^{-1}(r)|^p\rt)\\
\le& C_p\lt(\lt(|\wt{X}_{r}^{0}
-\wt{X}_{r}^{1}|^2+|t_0-t_1|^2\rt)|M_r|^{p-2}+|\wt{X}_{r}^{0}
-\wt{X}_{r}^{1}|^p+|t_0-t_1|^p\rt).
\end{aligned}
$$
It follows that, for $s\in[t,T]$,
\begin{eqnarray*}
&&\E\lt[|M_s|^p|\F_t\rt]\nonumber\\&=&\E\lt[\int_t^s\int_E
\lt(|M_r+\Delta \beta(r,e)|^p-|M_r|^p-p|M_r|^{p-2}M_r\Delta \beta(r,e)\rt)\Pi(\d e)\d r\big|\F_t\rt]\nonumber\\
&\le&C_p\E\lt[\int_t^s\lt(\lt(|\wt{X}_{r}^{0}
-\wt{X}_{r}^{1}|^2+|t_0-t_1|^2\rt)|M_r|^{p-2}+|\wt{X}_{r}^{0}
-\wt{X}_{r}^{1}|^p+|t_0-t_1|^p\rt)\d r\Big|\F_t\rt]\\
&\le&C_p\E\lt[\int_t^s\lt(|\wt{X}_{r}^{0}
-\wt{X}_{r}^{1}|^p+|t_0-t_1|^p\rt)\d r\Big|\F_t\rt]+C_p\E\lt[\int^s_t|M_r|^{p}\d r\Big|\F_t\rt],\nonumber
\end{eqnarray*}
and from Gronwall's inequality we obtain
\begin{equation}\label{eq:M}
\E\lt[|M_s|^p|\F_t\rt]\le C_p \E\lt[\int_t^s\lt(|\wt{X}_{r}^{0}
-\wt{X}_{r}^{1}|^p+|t_0-t_1|^p\rt)\d r\Big|\F_t\rt], \ s\in[t,T].
\end{equation}
Noticing that, for any
$t_\lambda\le t\le v\le T$,
\begin{eqnarray}\label{eq:x0-x1}
&&\wt{X}_v^{{0}}
-\wt{X}_v^{{1}}\nonumber\\
&=&\wt{X}_t^{{0}} -\wt{X}_t^{{1}}+
\int^v_{t}\lt(\frac{1}{\dot{\tau_0}}
b(\tau^{-1}_0(r),\widetilde{X}_r^{0},u^\lambda_r)-
\frac{1}{\dot{\tau_1}}b(\tau^{-1}_1(r),\widetilde{X}_r^{1},u^\lambda_r)\rt)\d
r \nonumber\\
&&+\int^v_{t}\lt(\frac{1}{\sqrt{\dot{\tau_0}}}
\sigma(\tau^{-1}_0(r),\widetilde{X}_r^{0},u^\lambda_r)-
\frac{1}{\sqrt{\dot{\tau_1}}}\sigma(\tau^{-1}_1(r),
\widetilde{X}_r^{1},u^\lambda_r)\rt)\d
B_r +M_v\\
&&+\int^v_{t}\hspace{-1mm}\int_E \lt(\lt(1-\frac{1}{\dot{\tau_0}}\rt)
\beta(\tau^{-1}_0(r),\widetilde{X}_{r}^{0},u^\lambda_{r},e)
-\lt(1-\frac{1}{\dot{\tau_1}}\rt)
\beta(\tau^{-1}_1(r),\widetilde{X}_r^{1},u^\lambda_{r},e)\rt)
\Pi(\d e)\d r,\nonumber
\end{eqnarray}
we get, by a standard argument (Recall that $\Pi(E)<\infty$) and Lemma \ref{lemma_tau}, the existence of a constant $C_{\delta,p}$ such that
\begin{eqnarray*}
&&\E\lt[\sup_{t\le v\le s}\lt|\wt{X}_v^{{0}}
-\wt{X}_v^{{1}}\rt|^p|\F_t\rt]\nonumber\\
&\le&C_{\delta,p}\lt|\wt{X}_t^{{0}}
-\wt{X}_t^{{1}}\rt|^p+C_p\E\lt[\lt|\int^s_{t}
\lt(\lt|\frac{1}{\dot{\tau}_0}-
\frac{1}{\dot{\tau}_1}\rt|+\lt|\frac{1}{\dot{\tau}_1}\rt|
\lt|\wt{X}_r^{{0}} -\wt{X}_r^{{1}}\rt|\rt)\d
r\rt|^p\bigg|\F_t\rt]\nonumber\\
&&\quad+C_{\delta,p}\E\lt[\lt(\int^s_{t} \lt(\lt|\frac{1}{\sqrt{\dot{\tau}_0}}-
\frac{1}{\sqrt{\dot{\tau}_1}}\rt|+\lt|\frac{1}{\sqrt{\dot{\tau}_1}}\rt|
\lt|\wt{X}_r^{{0}} - \wt{X}_r^{{1}}\rt|\rt)^2\d
r\rt)^{p/2}\bigg|\F_t\rt]\nonumber\\
&&\quad+C_{\delta,p}\E\lt[\lt|M_s\rt|^p|\F_t\rt]
%&&\quad+C\E\lt[\int^s_{t}\int_E\lt|\frac{1}{\dot{\tau}_0}
%-\frac{1}{\dot{\tau}_1}\rt|^p
%\lt|\beta(\tau^{-1}_0(r),\widetilde{X}_r^{0},
%u^\lambda_r,e)\rt|^p\Pi(\d
%e)\d r\bigg|\F_t\rt]\nonumber\\
+C_{\delta,p}|t_0-t_1|^p.
%\E\lt[\int^s_{t}\int_E\lt|1-\frac{1}{\dot{\tau}}\rt|^p
%\lt|\Delta\beta(r,e)\rt|^p
%\Pi(\d e)\d r\bigg|\F_t\rt].
%&\le&C\lt(\lt|\wt{X}_t^{{0}}
%-\wt{X}_t^{{1}}\rt|^p+|t_0-t_1|^p\rt)+C\int^s_{t}
%\E\lt[\lt|\wt{X}_r^{{0}}
%-\wt{X}_r^{{1}}\rt|^p\Big|\F_t\rt]\d r\nonumber\\
%&&+C\E\lt[\lt|\int^s_{t}\int_E\lt(\beta(\tau^{-1}_0(r),
%\widetilde{X}_{r-}^{0},u^\lambda_r,e)-
%\beta(\tau^{-1}_1(r),\widetilde{X}_{r-}^{1},
%u^\lambda_r,e)\rt)\tilde{\mu}(\d
%r,\d e)\rt|^p\bigg|\F_t\rt].
\end{eqnarray*}
Thus, taking into account (\ref{eq:M}), we obtain
\begin{eqnarray*}
&&\E\lt[\sup_{t\le v\le s}\lt|\wt{X}_v^{0}
-\wt{X}_v^{1}\rt|^p\Big|\F_t\rt]\\
&\le& C_{\delta,p}\lt(\lt|\wt{X}_t^{0}
-\wt{X}_t^{1}\rt|^p+|t_0-t_1|^p\rt)+C_{\delta,p}\int^s_t\E\lt[\sup_{t\le v\le r}\lt|\wt{X}_v^{0}
-\wt{X}_v^{1}\rt|^p\big|\F_t\rt]\d r,\ s\in[t_\lambda,T],
\end{eqnarray*}
and, finally, Gronwall's lemma yields that
\begin{eqnarray*}
\E\lt[\sup_{t\le v\le T}\lt|\wt{X}_v^{0}
-\wt{X}_v^{1}\rt|^p\Big|\F_t\rt]\le C_{\delta,p}\lt(|\wt{X}_t^{0}
-\wt{X}_t^{1}|^p+|t_0-t_1|^p\rt).\nonumber
\end{eqnarray*}
The proof is complete now. \hfill$\cajita$
\bigskip

Let us now prove Lemma \ref{lemma_lam-lam}:

\noindent\textit{Proof of Lemma \ref{lemma_lam-lam}:}
We observe that, for $s\in[t_\lambda,T],$
\begin{equation}\label{eq:x-lam-lam}
\begin{aligned}
&\wt{X}_s^{\lambda}-{X}_s^{{\lambda}}\\
= &\int^s_{t_\lambda}\lt(\frac{\lambda}{\dot{\tau_0}}
b(\tau^{-1}_0(r),\widetilde{X}_r^0,u^\lambda_r)+
\frac{1-\lambda}{\dot{\tau_1}}b(\tau^{-1}_1(r),\widetilde{X}_r^1,u^\lambda_r)
-b(r,{X}_r^{\lambda},u^\lambda_r)\rt)\d r\\
&+\int^s_{t_\lambda}\lt(\frac{\lambda}{\sqrt{\dot{\tau_0}}}
\sigma(\tau^{-1}_0(r),\widetilde{X}_r^0,u^\lambda_r)+
\frac{1-\lambda}{\sqrt{\dot{\tau_1}}}\sigma(\tau^{-1}_1(r),
\widetilde{X}_r^1,u^\lambda_r)
-\sigma(r,{X}_r^{\lambda},u^\lambda_r)\rt)\d B_r\\
&+\int^s_{t_\lambda}\int_E\lt(\lambda\beta(\tau^{-1}_0(r),
\widetilde{X}_{r-}^0,u^\lambda_r,e)+
(1-\lambda)\beta(\tau^{-1}_1(r),\widetilde{X}_{r-}^1,u^\lambda_r,e)-
\beta(r,{X}_{r-}^{\lambda},u^\lambda_r,e)\rt)\tilde{\mu}(\d
r,\d e)\\
&+\int^s_{t_\lambda}\hspace{-1mm}\int_E\lt(\lambda\lt(1-\frac{1}{\dot{\tau_0}}\rt)
\beta(\tau^{-1}_0(r),\widetilde{X}_r^0,u^\lambda_r,e)
+(1-\lambda)\lt(1-\frac{1}{\dot{\tau_1}}\rt)
\beta(\tau^{-1}_1(r),\widetilde{X}_r^1,u^\lambda_r,e)\rt)
\Pi(\d e)\d r.
\end{aligned}
\end{equation}
Taking into account that, as a consequence of the assumptions on the coefficients,
we have on one hand that the functions $b,\sigma,\beta$ but also $-b,-\sigma,-\beta$ are semiconcave in $(t,x)$,
uniformly with respect to $u$ and $(u,e)$, respectively, and that, on the other hand, $
\lambda(1-\frac{1}{\dot{\tau}_0})=-(1-\lambda)(1-\frac{1}{\dot{\tau}_1})
=\frac{\lambda(1-\lambda)}{T-t_\lambda}(t_0-t_1)
$ (see Lemma \ref{lemma_tau}), we deduce
\begin{eqnarray*}
&&\lt|\frac{\lambda}{\dot{\tau_0}}
b(\tau^{-1}_0(r),\widetilde{X}_r^0,u^\lambda_r)+
\frac{1-\lambda}{\dot{\tau_1}}b(\tau^{-1}_1(r),\widetilde{X}_r^1,u^\lambda_r)
-b(r,{X}_r^{\lambda},u^\lambda_r)\rt|\nonumber\\
&\le&\lt|{\lambda}
b(\tau^{-1}_0(r),\widetilde{X}_r^0,u^\lambda_r)+
{(1-\lambda)}b(\tau^{-1}_1(r),\widetilde{X}_r^1,u^\lambda_r)
-b(r,{X}_r^{\lambda},u^\lambda_r)\rt|\nonumber\\
&&+\lt|\lambda\lt(\frac{1}{\dot{\tau_0}}-1\rt)
b(\tau^{-1}_0(r),\widetilde{X}_r^0,u^\lambda_r)+
(1-\lambda)\lt(\frac{1}{\dot{\tau_1}}-1\rt)b(\tau^{-1}_1(r),\widetilde{X}_r^1,u^\lambda_r)\rt|
\nonumber\\
&\le&C_\delta\lambda(1-\lambda)\lt(|t_0-t_1|^2+\lt|\wt{X}_r^{0}
-\wt{X}_r^{1}\rt|^2\rt)
+C_\delta\lt|\wt{X}_r^{\lambda}-{X}_r^{\lambda}\rt|, \ r\in[t_\lambda,T].\nonumber
\end{eqnarray*}
Similarly, we get
\begin{eqnarray*}
&&\lt|\frac{\lambda}{\sqrt{\dot{\tau_0}}}
\sigma(\tau^{-1}_0(r),\widetilde{X}_r^0,u^\lambda_r)+
\frac{1-\lambda}{\sqrt{\dot{\tau_1}}}\sigma(\tau^{-1}_1(r),\widetilde{X}_r^1,u^\lambda_r)
-\sigma(r,{X}_r^{\lambda},u^\lambda_r)\rt|\nonumber\\
&\le&C_\delta\lambda(1-\lambda)\lt(|t_0-t_1|^2+\lt|\wt{X}_r^{0}
-\wt{X}_r^{1}\rt|^2\rt)
+C_\delta\lt|\wt{X}_r^{\lambda}-{X}_r^{\lambda}\rt|,\ r\in[t_\lambda,T],\nonumber
\end{eqnarray*}
and
\begin{eqnarray*}
&&\lt|{\lambda}
\beta(\tau^{-1}_0(r),\widetilde{X}_r^0,u^\lambda_r,e)+
(1-\lambda)\beta(\tau^{-1}_1(r),\widetilde{X}_r^1,u^\lambda_r,e)
-\beta(r,{X}_r^{\lambda},u^\lambda_r,e)\rt|\nonumber\\
&\le&C_\delta\lambda(1-\lambda)\lt(|t_0-t_1|^2+\lt|\wt{X}_r^{0}
-\wt{X}_r^{1}\rt|^2\rt)
+C_\delta\lt|\wt{X}_r^{\lambda}-{X}_r^{\lambda}\rt|,\ r\in[t_\lambda,T].\nonumber
\end{eqnarray*}
Moreover, by using again that $
\lambda(1-\frac{1}{\dot{\tau}_0})=-(1-\lambda)(1-\frac{1}{\dot{\tau}_1})
=\frac{\lambda(1-\lambda)}{T-t_\lambda}(t_0-t_1)
$, we obtain 
$$
\begin{aligned}
&\lt|\lambda\lt(1-\frac{1}{\dot{\tau_0}}\rt)
\beta(\tau^{-1}_0(r),\widetilde{X}_r^0,u^\lambda_r,e)
+(1-\lambda)\lt(1-\frac{1}{\dot{\tau_1}}\rt)
\beta(\tau^{-1}_1(r),\widetilde{X}_r^1,u^\lambda_r,e)\rt|\\
\le&
C_\delta\lambda(1-\lambda)\lt(|t_0-t_1|^2+\lt|\wt{X}_r^{0}
-\wt{X}_r^{1}\rt|^2\rt), \ r\in[t_\lambda,T].
\end{aligned}
$$
%Taking into account of
%$
%\lambda(1-\frac{1}{\dot{\tau}_0})=-(1-\lambda)(1-\frac{1}{\dot{\tau}_1})
%=\frac{\lambda(1-\lambda)}{T-t_\lambda}(t_0-t_1),
%$
%and applying the same method as that used in the proof of  Lemma \ref{lemma_1-2},  we obtain
Consequently, by combining the above estimates with the argument which has lead to (\ref{eq:M}) 
in the proof of Lemma \ref{lemma_1-2}, we get 
\begin{eqnarray*}
&&\E\lt[\sup_{t\le v\le s}\lt|\wt{X}_v^{\lambda}-{X}_v^{{\lambda}}
\rt|^p\bigg|\F_t\rt]\le
\nonumber\\
&\le&C_{\delta,p}\lt|\wt{X}_t^{\lambda}-{X}_t^{{\lambda}}
\rt|^p+C_\delta(\lambda(1-\lambda))^p\lt(\lt|t_0-t_1\rt|^{2p}+\E\lt[\int^s_{t}
\lt|\widetilde{X}_r^0-\widetilde{X}_r^1\rt|^{2p}\d
r\bigg|\F_t\rt]\rt)\\
&&+C_{\delta,p}\int^s_{t}\E\lt[\sup_{t\le v\le r}
\lt|\wt{X}_v^{\lambda}-{X}_v^{\lambda}\rt|^p\bigg|\F_t\rt]\d r,\ \ t_\lambda\le t\le s\le T,
\end{eqnarray*}
and, thus, Lemma \ref{lemma_1-2} yields 
\begin{eqnarray*}
\E\lt[\sup_{t\le v\le s}\lt|\wt{X}_v^{\lambda}-{X}_v^{{\lambda}}
\rt|^p\bigg|\F_t\rt]&\le&C_{\delta,p}\lt|\wt{X}_t^{\lambda}-{X}_t^{{\lambda}}
\rt|^p+C_{\delta,p}(\lambda(1-\lambda))^p\lt(|t_0-t_1|^{2p}+\lt|\wt{X}_t^{0}
-\wt{X}_t^{1}\rt|^{2p}\rt)\\
&&+C_{\delta,p}\int^s_{t}\E\lt[\sup_{t\le v\le r}\lt|\wt{X}_v^{\lambda}-{X}_v^{\lambda}\rt|^p\bigg|\F_t\rt]\d
r,\ \ t_\lambda\le t\le s\le T.
\end{eqnarray*}
Hence, Gronwall's inequality gives
\[
\E\lt[\sup_{t\le s\le T}\lt|\wt{X}_s^{\lambda}-{X}_s^{\lambda}\rt|^p\bigg|\F_t\rt]\le
C_{\delta,p}\lt|\wt{X}_t^{\lambda}-{X}_t^{{\lambda}}
\rt|^p+C_{\delta,p}(\lambda(1-\lambda))^p\lt(|t_0-t_1|^{2p}+\lt|\wt{X}_t^{0}
-\wt{X}_t^{1}\rt|^{2p}\rt).
\]
%Hence, we have
%\begin{eqnarray*}
%&&\E\lt[\sup_{t\le s\le T}
%\lt|\wt{X}_s^{\lambda}-{X}_s^{\lambda}\rt|^p\bigg|\F_t\rt]\nonumber\\
%&\le&
%\E\lt[\sup_{t\le s
%\le
%T}\lt|\wt{X}_s^{\lambda}-{X}_s^{\lambda}\rt|^{2p}\bigg|\F_t\rt]^{1/2}\nonumber\\
%&\le&
%C_{\delta,p}\lt|\wt{X}_t^{\lambda}-{X}_t^{{\lambda}}
%\rt|^p+C_{\delta,p}(\lambda(1-\lambda))^p\lt(|t_0-t_1|^{2p}+\lt|\wt{X}_t^{0}
%-\wt{X}_t^{1}\rt|^{2p}\rt), \ t\in[t_\lambda,T].\nonumber
%\end{eqnarray*}
The proof of Lemma \ref{lemma_lam-lam} is complete.\hfill$\cajita$
\bigskip

We prove now the analogue estimates for our BSDEs stated in Lemma \ref{lemma_Y0-Y1}.

\noindent\textit{Proof of Lemma \ref{lemma_Y0-Y1}:}
We notice that, for $t_\lambda\le s\le t\le T$,
\begin{equation}\label{eq:f}
\begin{aligned}
&\wt{Y}^{0}_s-\wt{Y}^1_s\\
=&\Phi(\wt{X}_T^{0})-\Phi(\wt{X}_T^{1}) \\
&+\int^T_s\lt(
\frac{1}{\dot{\tau}_0}f(\tau^{-1}_0(r),\wt{X}_r^{0},
\wt{Y}_r^{0},\sqrt{\dot{\tau}_0}\wt{Z}_r^{0},
\wt{U}_r^{0},u^\lambda_r)-\frac{1}{\dot{\tau}_1}f(\tau^{-1}_1(r),\wt{X}_r^{1},
\wt{Y}_r^{1},\sqrt{\dot{\tau}_1}\wt{Z}_r^{1},
\wt{U}_r^{1},u^\lambda_r)\rt)\d
r\\&-\int^T_s \lt(\wt{Z}_r^{0}-\wt{Z}_r^{1}\rt)\d B_r-
\int^T_s\int_E\lt(\wt{U}_r^{0}(e)-\wt{U}_r^{1}(e)\rt)\widetilde{\mu}(\d
r, \d e)\\
&- \int^T_s\int_E\lt(\lt(1-\frac{1}{\dot{\tau_0}}\rt)\wt{U}_r^{0}(e)
-\lt(1-\frac{1}{\dot{\tau_1}}\rt)\wt{U}_r^{1}(e)\rt)\Pi (\d e)\d r.
\end{aligned}
\end{equation}
By applying It\^o's formula to $|\wt{Y}^0_s-\wt{Y}^1_s|^2$, we get
from the boundedness and the Lipschitz continuity of $f$ that,
for $t_\lambda \le s\le t\le T$,
\begin{eqnarray}\label{eq:y0-y1}
&&|\wt{Y}^0_s-\wt{Y}^1_s|^2
+\int^T_s\lt|\wt{Z}_r^{0}-\wt{Z}_r^{1}\rt|^2\d r
+\int^T_s\int_E\lt|\wt{U}_r^{0}(e)-\wt{U}_r^{1}(e)\rt|^2\Pi(\d e) \d r\nonumber\\
&\le&\lt|\Phi(\wt{X}_T^{0})-\Phi(\wt{X}_T^{1})\rt|^2
+C|t_0-t_1|^2+C\int^T_s\lt|\wt{X}_r^{0}-\wt{X}_r^{1}\rt|^2\d r
+C\int^T_s|\wt{Y}^0_r-\wt{Y}^1_r|^2\d r\nonumber\\
&&-2\int^T_s\lt(\wt{Y}^0_r-\wt{Y}^1_r\rt)\lt(\wt{Z}^0_r-\wt{Z}^1_r\rt)\d B_r
\nonumber\\
&&+|t_0-t_1|^2\int^T_s\lt(\lt|\wt{Z}_r^1\rt|^2
+\int_E\lt(\lt|\wt{U}_r^{0}(e)\rt|^2+\lt|\wt{U}_r^{1}(e)\rt|^2\rt)\Pi(\d e)\rt)\d r\nonumber\\
&&-\int^T_s\int_E\lt(2\lt(\wt{Y}^0_r-\wt{Y}^1_r\rt)\lt(\wt{U}^0_r(e)-\wt{U}^1_r(e)\rt)
+ \lt|\wt{U}_r^{0}(e)-\wt{U}_r^{1}(e)\rt|^2\rt)\tilde{\mu}(\d r,\d e).
\end{eqnarray}
Then  taking the conditional expectation and applying Gronwall's inequality
and the Burkholder-Davis-Gundy inequality we obtain from the Lemmata \ref{lemma_apriori}
and \ref{lemma_1-2} (Recall also the arguments given in the proof of Lemma \ref{lemma_lam-lam}) that
\begin{eqnarray*}
&&\E\lt[\sup_{t\le r\le T}|\wt{Y}^0_r-\wt{Y}^1_r|^2
+\int^T_t\lt|\wt{Z}_r^{0}-\wt{Z}_r^{1}\rt|^2\d r
+\int^T_t\int_E\lt|\wt{U}_r^{0}(e)-\wt{U}_r^{1}(e)\rt|^2\Pi(\d e) \d r\bigg|\F_t\rt]\nonumber\\
&\le& C_\delta\lt(\lt|\wt{X}_t^{0}-\wt{X}_t^{1}\rt|^2+|t_0-t_1|^2\rt),\ t\in[t_\lambda,T].\nonumber
\end{eqnarray*}
In particular, we have
\[
|\wt{Y}^0_t-\wt{Y}^1_t|^2\le
C_\delta\lt(\lt|\wt{X}_t^{0}-\wt{X}_t^{1}\rt|^2+|t_0-t_1|^2\rt), \ \ t\in[t_\lambda,T].
\]
Consequently, from Lemma \ref{lemma_lam-lam}, we have, for $p\ge 2$,
\begin{eqnarray}\label{eq:supY}
\E\lt[\sup_{t\le r\le T}|\wt{Y}^0_r-\wt{Y}^1_r|^p\bigg|\F_t\rt]
\le C_\delta\lt(\lt|\wt{X}_t^{0}-\wt{X}_t^{1}\rt|^p+|t_0-t_1|^p\rt),\ t\in[t_\lambda,T].
\end{eqnarray}
From (\ref{eq:y0-y1}), by using that $\tilde{\mu}(\d r \d e)=\mu(\d r\d e)-\Pi(\d e)\d r$, we get also that
\begin{eqnarray*}
&&\int^T_s\lt|\wt{Z}_r^{0}-\wt{Z}_r^{1}\rt|^2\d r
+\int^T_s\int_E\lt|\wt{U}_r^{0}(e)-\wt{U}_r^{1}(e)\rt|^2\mu(\d r,\d e)\nonumber\\
&\le&\lt|\Phi(\wt{X}_T^{0})-\Phi(\wt{X}_T^{1})\rt|^2
+C|t_0-t_1|^2+C\int^T_s\lt|\wt{X}_r^{0}-\wt{X}_r^{1}\rt|^2\d r
+C\int^T_s|\wt{Y}^0_r-\wt{Y}^1_r|^2\d r\nonumber\\
&&-2\int^T_s\lt(\wt{Y}^0_r-\wt{Y}^1_r\rt)\lt(\wt{Z}^0_r-\wt{Z}^1_r\rt)\d B_r
+|t_0-t_1|^2\int^T_s\lt(\lt|\wt{Z}_r^1\rt|^2
+\int_E\lt(\lt|\wt{U}_r^{0}(e)\rt|^2+\lt|\wt{U}_r^{1}(e)\rt|^2\rt)\Pi(\d e)\rt)\d r\nonumber\\
&&-\int^T_s\int_E\lt(2\lt(\wt{Y}^0_r-\wt{Y}^1_r\rt)\lt(\wt{U}^0_r(e)-\wt{U}^1_r(e)\rt)\rt)\tilde{\mu}(\d r,\d e).
\end{eqnarray*}
Hence, by applying
Lemma \ref{lemma_1-2} and Burkholder-Davis-Gundy's
inequality, we obtain from the above inequality that
\begin{eqnarray}\label{eq:supZU}
&&\E\lt[\lt(\int^T_s\lt|\wt{Z}_r^{0}-\wt{Z}_r^{1}\rt|^2\d r\rt)^{{p}/{2}}
+\lt(\int^T_s\int_E\lt|\wt{U}_r^{0}(e)-\wt{U}_r^{1}(e)\rt|^2\mu(\d r,\d e)\rt)
^{{p}/{2}}\bigg|\F_t\rt]\nonumber\\
&\le& C\lt(\lt|\wt{X}_t^{0}-\wt{X}_t^{1}\rt|^p+|t_0-t_1|^p\rt)
+C\E\bigg[\lt(\int^T_s\lt|\wt{Y}^0_r-\wt{Y}^1_r\rt|^2
\lt|\wt{Z}^0_r-\wt{Z}^1_r\rt|^2\d r\rt)^{{p}/{4}}\nonumber\\
&&\qquad+2\lt(\int^T_s\int_E\lt|\wt{Y}^0_r-\wt{Y}^1_r\rt|^2
\lt|\wt{U}^0_r(e)-\wt{U}^1_r(e)\rt|^2{\mu}(\d r,\d e)\rt)^{{p}/{4}}\bigg|\F_t\bigg]\nonumber\\
&\le&\hspace{-1mm}C\lt(\lt|\wt{X}_t^{0}-\wt{X}_t^{1}\rt|^p+|t_0-t_1|^p\rt)
+\E\bigg[C\sup_{s\le r\le T}\lt|\wt{Y}^0_r-\wt{Y}^1_r\rt|^p
+\frac{1}{2}\lt(\int^T_s\lt|\wt{Z}^0_r-\wt{Z}^1_r\rt|^2\d r\rt)^{{p}/{2}}\nonumber\\
&&\qquad+\frac12\lt(\int^T_s\int_E\lt|\wt{U}^0_r(e)-\wt{U}^1_r(e)\rt|^2{\mu}
(\d r,\d e)\rt)^{{p}/{2}}\bigg|\F_t\bigg].
\end{eqnarray}
Combining (\ref{eq:supY}) and (\ref{eq:supZU}), we get
\begin{eqnarray}\label{eq:estYY}
&&\E\bigg[\sup_{s\le r\le T}|\wt{Y}^0_r-\wt{Y}^1_r|^p
+\lt(\int^T_s\lt|\wt{Z}^0_r-\wt{Z}^1_r\rt|^2\d r\rt)^{p/2} + \lt(\int^T_s\int_E\lt|\wt{U}^0_r(e)-\wt{U}^1_r(e)\rt|^2{\mu}
(\d r,\d e)\rt)^{p/2}\bigg|\F_t\bigg]\nonumber\\
&\le& C\lt(|\wt{X}_t^{0}
-\wt{X}_t^{1}|^p+|t_0-t_1|^p\rt).
\end{eqnarray}
In order to replace in the left-hand side of (\ref{eq:estYY}) the term 
$\E\bigg[\lt(\int^T_s\int_E\lt|\wt{U}^0_r(e)-\wt{U}^1_r(e)\rt|^2{\mu}
(\d r,\d e)\rt)^{p/2}\bigg|\F_t\bigg]$ 
by 
$\E\bigg[\lt(\int^T_s\int_E\lt|\wt{U}^0_r(e)-\wt{U}^1_r(e)\rt|^2{\Pi}
(\d e)\d r\rt)^{p/2}\bigg|\F_t\bigg]$, we make the following estimate
\begin{eqnarray*}
&&\E\lt[\lt(\int^T_s\int_E\lt|\wt{U}^0_r(e)-\wt{U}^1_r(e)\rt|^2
\Pi(\d e) \d r\rt)^{p/2}\bigg|\F_t\rt]\nonumber\\
&\le& (\Pi(E)T)^{\frac{p-2}{2}}
\E\lt[\int^T_s\int_E\lt|\wt{U}^0_r(e)-\wt{U}^1_r(e)\rt|^p
\Pi(\d e)\d r\bigg|\F_t\rt]\nonumber\\
&=& (\Pi(E)T)^{\frac{p-2}{2}}
\E\lt[\int^T_s\int_E\lt|\wt{U}^0_r(e)-\wt{U}^1_r(e)\rt|^p
{\mu}(\d r,\d e)\bigg|\F_t\rt].
\end{eqnarray*}
Let us denote by $N$ the Poisson process $N_s=\mu([t_\lambda,s]\times E), s\in[t_\lambda,T]$, with intensity $\Pi(E)$, the associated sequence of jump times by $\tau_i=\inf\{s\ge t_\lambda: N_s=i\}, i\ge 1$, and by $p_i: \{\tau_i\le T\}\to E$ s.t. $\mu(\{(\tau_i,p_i)\})=1$ on $\{\tau_i\le T\}$, the associated sequence of marks. We observe that
\begin{eqnarray*}
&&\int^T_s\int_E\lt|\wt{U}^0_r(e)-\wt{U}^1_r(e)\rt|^p
{\mu}(\d r,\d e)
= \sum_{i\ge 1}
\lt|\wt{U}^0_{\tau_i}(p_i)-\wt{U}_{\tau_i}^{1}(p_i)\rt|^p\\
&\le& \lt(\sum_{i\ge 1} \lt|\wt{U}^0_{\tau_i}(p_i)-\wt{U}_{\tau_i}^{1}(p_i)\rt|^2\rt)^{p/2}
=\lt(\int^T_s\int_E\lt|\wt{U}^0_r(e)-\wt{U}^1_r(e)\rt|^2
{\mu}(\d r,\d e)\rt)^{p/2}.
\end{eqnarray*}
Consequently, from the previous estimate we get
\begin{eqnarray*}
&&\E\lt[\lt(\int^T_s\int_E\lt|\wt{U}^0_r(e)-\wt{U}^1_r(e)\rt|^2
\Pi(\d e) \d r\rt)^{p/2}\bigg|\F_t\rt]\nonumber\\
%&\le& (\Pi(E)T)^{\frac{p-2}{2}}
%\E\lt[\int^T_s\int_E\lt|\wt{U}^0_r(e)-\wt{U}^1_r(e)\rt|^p
%\Pi(\d e)\d r\bigg|\F_t\rt]\nonumber\\
&\le& (\Pi(E)T)^{\frac{p-2}{2}}
\E\lt[\int^T_s\int_E\lt|\wt{U}^0_r(e)-\wt{U}^1_r(e)\rt|^p
{\mu}(\d r,\d e)\bigg|\F_t\rt]\\
%&=& (\Pi(E)T)^{\frac{p-2}{2}}\E\lt[\sum_{s\le r\le T}
%\lt|\wt{U}^0_r(P_r)-\wt{U}_r^{1}(P_r)\rt|^p\bigg|\F_t\rt]\\
%&\le& (\Pi(E)T)^{\frac{p-2}{2}}\E\lt[\lt(\sum_{s\le r\le T} \lt|\wt{U}^0_r(P_r)
%-\wt{U}_r^{1}(P_r)\rt|^2\rt)^{p/2}\bigg|\F_t\rt]
%\\
%\E\bigg[\int^T_s\int_E\lt(\int^T_s\int_E
%\lt|\wt{U}^0_r(e)-\wt{U}^1_r(e)\rt|^2{\mu}(\d r,\d e)\rt)^{\frac{p}{2}-1}
%\nonumber\\
%&&\qquad\qquad\qquad
%\times\lt|\wt{U}^0_r(e)-\wt{U}^1_r(e)\rt|^2{\mu}(\d r,\d e)\bigg|\F_t\bigg]
%\nonumber\\
&\le&(\Pi(E)T)^{\frac{p-2}{2}}\E\lt[\lt(\int^T_s\int_E
\lt|\wt{U}^0_r(e)-\wt{U}^1_r(e)\rt|^2
{\mu}(\d r,\d e)\rt)^{p/2}\bigg|\F_t\rt].
\end{eqnarray*}
This latter estimate allows to deduce from (\ref{eq:estYY}) the wished result.\hfill$\cajita$
\bigskip

Let us now give the proof of Lemma \ref{lemma_comparison}.

\noindent\textit{Proof of Lemma \ref{lemma_comparison}}:
%From the semiconcavity of $f$ and Lemma \ref{lemma_tau} we obtain,
First we notice that,
for $r\in[t_\lambda,T]$,
\begin{equation*}
\begin{aligned}
I:=&\frac{\lambda}{\dot{\tau_0}}f(\tau_0^{-1}(r),\wt{X}_r^0,
\wt{Y}_r^0,\sqrt{\dot{\tau_0}}\wt{Z}_r^0,\wt{U}_r^0,u^\lambda_r)
+\frac{1-\lambda}{\dot{\tau_1}}f(\tau_1^{-1}(r),\wt{X}_r^1,
\wt{Y}_r^1,\sqrt{\dot{\tau_1}}\wt{Z}_r^1,\wt{U}_r^1,u^\lambda_r)\\
=&\lambda f(\tau_0^{-1}(r),\wt{X}_r^0,\wt{Y}_r^0,\sqrt{\dot{\tau_0}}
\wt{Z}_r^0,\wt{U}_r^0,u^\lambda_r)+(1-\lambda) f(\tau_1^{-1}(r),
\wt{X}_r^1,\wt{Y}_r^1,\sqrt{\dot{\tau_1}}\wt{Z}_r^1,\wt{U}_r^1,u^\lambda_r)\\
&-\lambda(1-\lambda)\frac{t_1-t_0}{T-t_\lambda}\Big(f(\tau_0^{-1}(r),
\wt{X}_r^0,\wt{Y}_r^0,\sqrt{\dot{\tau_0}}\wt{Z}_r^0,\wt{U}_r^0,u^\lambda_r)
-f(\tau_1^{-1}(r),\wt{X}_r^1,\wt{Y}_r^1,\sqrt{\dot{\tau_1}}\wt{Z}_r^1,
\wt{U}_r^1,u^\lambda_r)\Big).
\end{aligned}
\end{equation*}
We also observe that, by applying Lemma \ref{lemma_tau},
$$
\begin{aligned}
&\lt|\lambda\sqrt{\dot{\tau_0}}
\wt{Z}_r^0+(1-\lambda)\sqrt{\dot{\tau_1}}\wt{Z}_r^1-\wt{Z}_r^\lambda\rt|\\
=&\lt|\lambda(\sqrt{\dot{\tau_0}}-1)
\wt{Z}_r^0+(1-\lambda)(\sqrt{\dot{\tau_1}}-1)\wt{Z}_r^1\rt|\\
\le & \lambda|1-\sqrt{\dot{\tau_0}}||\wt{Z}_r^0-\wt{Z}_r^1|+|\lambda(1-\sqrt{\dot{\tau_0}})
+(1-\lambda)(1-\sqrt{\dot{\tau_1}})||\wt{Z}_r^1|\\
\le& C_\delta\lambda(1-\lambda)(|t_0-t_1||\wt{Z}_r^0-\wt{Z}_r^1|+|t_0-t_1|^2|\wt{Z}_r^1|).
\end{aligned}
$$
Thus, by using the estimates from  Lemma \ref{lemma_tau}, 
the semiconcavity as well as the Lipschitz continuity of $f$,
we have
\begin{eqnarray*}
I&\le&f(r,\wt{X}_r^\lambda,\wt{Y}_r^\lambda,\lambda\sqrt{\dot{\tau_0}}
\wt{Z}_r^0+(1-\lambda)\sqrt{\dot{\tau_1}}\wt{Z}_r^1,\wt{U}_r^\lambda,u^\lambda_r)
+C_\delta \lambda(1-\lambda)\bigg(|t_0-t_1|^2\lt(1+\lt|\wt{Z}_r^1\rt|^2\rt)\\
&&+\lt|\wt{X}_r^0-\wt{X}_r^1\rt|^2+\lt|\wt{Y}_r^0-\wt{Y}_r^1\rt|^2+
\lt|\wt{Z}_r^0-\wt{Z}_r^1\rt|^2+
\int_E \lt|\wt{U}_r^0(e)-\wt{U}_r^1(e)\rt|^2\Pi(\d e)\bigg)\\
&\le&f(r,\wt{X}_r^\lambda,\wt{Y}_r^\lambda,\wt{Z}_r^\lambda,
\wt{U}_r^\lambda,u^\lambda_r)+C_\delta \lambda(1-\lambda)
\bigg(|t_0-t_1|^2\lt(1+\lt|\wt{Z}_r^1\rt|^2\rt)+\lt|\wt{X}_r^0-\wt{X}_r^1\rt|^2\\
&&+
\lt|\wt{Y}_r^0-\wt{Y}_r^1\rt|^2+\lt|\wt{Z}_r^0-\wt{Z}_r^1\rt|^2+
\int_E \lt|\wt{U}_r^0(e)-\wt{U}_r^1(e)\rt|^2\Pi(\d e)\bigg).
\end{eqnarray*}
Finally, by taking into account the Lipschitz continuity of $f$
as well as the definition of $D_r$, we get
\begin{eqnarray*}
I&\le&f(r,\wt{X}_r^\lambda,\wt{Y}_r^\lambda-D_r,
\wt{Z}_r^\lambda,\wt{U}_r^\lambda,u^\lambda_r)+C_\delta \lambda(1-\lambda)
\bigg(|t_0-t_1|^2\lt(1+\lt|\wt{Z}_r^1\rt|^2\rt)
+\lt|\wt{X}_r^0-\wt{X}_r^1\rt|^2\nonumber\\
&&+
\lt|\wt{Y}_r^0-\wt{Y}_r^1\rt|^2+\lt|\wt{Z}_r^0-\wt{Z}_r^1\rt|^2+
\int_E \lt|\wt{U}_r^0(e)-\wt{U}_r^1(e)\rt|^2\Pi(\d e)\bigg)+CD_r.\nonumber
\end{eqnarray*}
Moreover, thanks to the semiconcavity of $\Phi$, we have
$$
\lambda\Phi\lt(\wt{X}_T^{0}\rt)+
(1-\lambda)\Phi\lt(\wt{X}_T^{1}\rt)\le
\Phi\lt(\wt{X}_T^{\lambda}\rt)+CB_T+C_\delta\lambda(1-\lambda)A_T^2.
$$
By virtue of  assumption $\textbf{(H5)}$, we can use
the comparison theorem in \cite{Ro} (Theorem 2.5) in order to conclude that
$$\wt{Y}_s^\lambda\le\wh{Y}_s^\lambda,\ \mathrm{for}\ s\in[t_\lambda,T].$$
The proof of Lemma \ref{lemma_comparison} is complete.\hfill$\cajita$
\bigskip

\noindent\textit{Proof of Lemma \ref{lemma_wtY-Y}:} For some
constant $\gamma>0$ which will be specified later, we apply Ito's
formula to $e^{\gamma s}(\overline{Y}_s^\lambda-Y^\lambda_s)^2$, and we get
\begin{eqnarray*}\label{eq:itoy}
&&\d e^{\gamma s} (\overline{Y}_s^\lambda-Y^\lambda_s)^2\\
&=&- 2e^{\gamma s} (\overline{Y}_s^\lambda-Y^\lambda_s)
\bigg[f(s,{X}_s^{\lambda},\overline{Y}_s^\lambda,\wh{Z}_s^\lambda,
\wh{U}_s^\lambda,u^\lambda_s)
-f(s,{X}_s^{\lambda},Y^\lambda_s,Z^\lambda_s,U^\lambda_s,u^\lambda_s)
+CD_s\\
&&+C_\delta^0\lambda(1-\lambda)\lt(|t_0-t_1|^2\lt(1+|\wt{Z}_s^1|^2\rt)
+|\wt{Z}_s^0-\wt{Z}_s^1|^2\rt)
+\int_E\lt|\wt{U}_s^{0}(e)-\wt{U}_s^{1}(e)\rt|^2\Pi(\d e)\bigg]
\d s\\
&&+2e^{\gamma s} (\overline{Y}_s^\lambda-Y^\lambda_s)
\lt(\wh{Z}_s^\lambda-Z^\lambda_s\rt)\d
B_s+2e^{\gamma s} (\overline{Y}_{s-}^\lambda-Y^\lambda_{s-})
\int_E\lt(\wh{U}_s^{\lambda}(e)-U_s^\lambda(e)\rt)\tilde{\mu}(\d
s, \d e)\\
&&-2e^{\gamma s} (\overline{Y}_{s-}^\lambda-Y^\lambda_{s-})\d D_s
+e^{\gamma s}\lt(\gamma\lt|\overline{Y}_s^\lambda-Y^\lambda_s\rt|^2
+\lt|\wh{Z}_s^\lambda-Z^\lambda_s\rt|^2\rt)\d s \\
&&+e^{\gamma s}\d\lt[\overline{Y}^\lambda-Y^\lambda\rt]_s^d, \ s\in[t_\lambda,T],
\end{eqnarray*}
where $\lt[\overline{Y}^\lambda-Y^\lambda\rt]_s^d$
denotes the purely discontinuous part of the
quadratic variation of $\overline{Y}^\lambda-Y^\lambda$:
$$
\lt[\overline{Y}^\lambda-Y^\lambda\rt]_s^d=
\sum_{t_\lambda<r\le s}\lt(\Delta\overline{Y}^\lambda_r
-\Delta Y^\lambda_r\rt)^2, \ s\in[t_\lambda,T].
$$

We notice that, for $s\in[t_\lambda,T]$,
\begin{eqnarray*}
&&\int^T_s e^{\gamma r}\d \lt[\overline{Y}^\lambda-Y^\lambda\rt]_r^d\\
&=&\sum_{s\le r\le T} e^{\gamma r}\lt(\Delta\overline{Y}^\lambda_r
-\Delta Y^\lambda_r\rt)^2\\
&=&\sum_{s\le r\le T} e^{\gamma r}\lt(\int_E\lt(\wh{U}_r^\lambda(e)-U^\lambda_r(e)\rt)
\mu(\{r\},\d e)-\Delta D_r\rt)^2\\
&\ge&\frac{1}{2}\int^T_s \int_E e^{\gamma r}\lt(\wh{U}_r^\lambda(e)-U^\lambda_r(e)\rt)^2
\mu(\d r,\d e)-C_{\gamma}\int^T_s|\Delta D_r|\d D_r\\
&\ge&\frac{1}{2}\int^T_s \int_E e^{\gamma r}\lt(\wh{U}_r^\lambda(e)-U^\lambda_r(e)\rt)^2
\mu(\d r,\d e)-C_{\gamma}(D_T-D_s)^2,
\end{eqnarray*}
where $C_{\gamma}=e^{\gamma T}$.
Hence, by integrating from $s\in[t_\lambda,T]$ to $T$ and taking
the conditional expectation on both sides, we deduce that
\begin{eqnarray*}
&& e^{\gamma s} \lt|\overline{Y}_s^\lambda-Y^\lambda_s\rt|^2+
\E\bigg[\int^T_se^{\gamma r}\lt(\gamma\lt|\overline{Y}_r^\lambda-Y^\lambda_r\rt|^2
+\lt|\wh{Z}_r^\lambda-Z^\lambda_r\rt|^2\rt)\d r\nonumber\\
&&\quad\qquad+\frac{1}{2}\int^T_s\int_E e^{\gamma r}
\lt|\wh{U}_r^{\lambda}(e)-U_r^\lambda(e)\rt|^2{\mu}(\d
r, \d e)\bigg| \F_s\bigg]\nonumber\\
&\le&\E\bigg[\int^T_s2e^{\gamma r} (\overline{Y}_r^\lambda-Y^\lambda_r)
\bigg\{f(r,{X}_r^{\lambda},
\overline{Y}_r^\lambda,\wh{Z}_r^\lambda,\wh{U}_r^\lambda,u^\lambda_r)
-f(r,{X}_r^{\lambda},Y^\lambda_r,Z^\lambda_r,U^\lambda_r,u^\lambda_r)
+CD_r\nonumber\\
&&\quad+C_\delta^0\lambda(1-\lambda)\lt(|t_0-t_1|^2\lt(1+|\wt{Z}_r^1|^2\rt)
+|\wt{Z}_r^0-\wt{Z}_r^1|^2\rt)\nonumber\\
&&
+\int_E\lt|\wt{U}_r^{0}(e)-\wt{U}_r^{1}(e)\rt|^2\Pi(\d e)\bigg\}
\d r+\int^T_s2e^{\gamma r} (\overline{Y}_{r-}^\lambda-Y^\lambda_{r-})\d D_r
+C_{\gamma}(D_T-D_s)^2\bigg|\F_s\bigg].
\end{eqnarray*}
Then, by a standard argument and from the Lipschitz continuity of $f$, we get
\begin{eqnarray*}
&& e^{\gamma s} \lt|\overline{Y}_s^\lambda-Y^\lambda_s\rt|^2+\E\bigg[\int^T_se^{\gamma r}\lt(\gamma\lt|\overline{Y}_r^\lambda-Y^\lambda_r\rt|^2+\lt|\wh{Z}_r^\lambda
-Z^\lambda_r\rt|^2\rt)\d r\nonumber\\
&&\quad\qquad+\frac{1}{2}\int^T_s\int_E e^{\gamma r}
\lt|\wh{U}_r^{\lambda}(e)-U_r^\lambda(e)\rt|^2{\mu}(\d
r, \d e)\bigg| \F_s\bigg]\nonumber\\
&\le&\E\bigg[C_K\int^T_se^{\gamma r} \lt|\overline{Y}_r^\lambda
-Y^\lambda_r\rt|^2\d r+\frac12\int^T_se^{\gamma r}
\lt|\wh{Z}_r^\lambda-Z^\lambda_r\rt|^2\d r\nonumber\\
&&+\frac14\int^T_s\int_E e^{\gamma r}
\lt|\wh{U}_r^{\lambda}(e)-U_r^\lambda(e)\rt|^2\Pi(\d e)\d r\bigg|\F_s\bigg]\nonumber\\
&&+C_{\delta,\gamma}\E\bigg[\sup_{s\le r\le T}
\lt|\overline{Y}_r^\lambda-Y^\lambda_r\rt|D_{s,T}
+C_\gamma (D_T-D_s)^2\bigg|\F_s\bigg],\nonumber
\end{eqnarray*}
where
$$
\begin{aligned}
D_{s,T}=&D_T+\lambda(1-\lambda)\lt(|t_0-t_1|^2
\lt(1+\int^T_s|\wt{Z}_r^1|^2\d r\rt)+\int^T_s|\wt{Z}_r^0-\wt{Z}_r^1|^2\d r\rt)\\
&+\int_s^T\int_E|\wt{U}_r^{0}(e)-\wt{U}_r^{1}(e)|^2\Pi(\d e)\d r.
\end{aligned}
$$
Therefore, by choosing $\gamma$ large enough and
applying Lemma \ref{lemma_1-2} and Corollary \ref{cor_D},
we have the following estimate:
\begin{eqnarray}\label{eq:_Y-Y^2}
&&\lt|\overline{Y}_s^\lambda-Y^\lambda_s\rt|^2
+\E\lt[\int^T_s\lt|\wh{Z}_r^\lambda-Z^\lambda_r\rt|^2\d r
+\int^T_s\int_E \lt|\wh{U}_r^{\lambda}(e)-U_r^\lambda(e)\rt|^2{\mu}(\d
r, \d e)\bigg| \F_s\rt]\nonumber\\
&\le&C_{\delta,\gamma}\E\bigg[\sup_{s\le r\le T}\lt|\overline{Y}_r^\lambda-Y^\lambda_r\rt|D_{s,T}
\bigg|\F_s\bigg]+C_\gamma D_s^2, \ s\in[t_\lambda,T].
\end{eqnarray}
From Lemma \ref{lemma_Y0-Y1} and Corollary \ref{cor_D},
we get that for $p\ge 2$,
\begin{equation}\label{eq:D_t}
\E\lt[|D_{t,T}|^p|\F_t\rt]\le C(|\wt{X}_t^0-\wt{X}_t^1|^{2p}+|t_0-t_1|^{2p})\le C D_t^p.
\end{equation}

Let $1<p<2$ and $q>2$ be two constants such that
$\frac{1}{p}+\frac{1}{q}=1$. Then we have, for $\ep>0$,
\begin{equation*}
\begin{aligned}
&\E\bigg[\sup_{s\le r\le T}\lt|\overline{Y}_r^\lambda-Y^\lambda_r\rt|D_{s,T}
\bigg|\F_s\bigg]\\
\le&\E\bigg[\sup_{s\le r\le T}\lt|\overline{Y}_r^\lambda-
Y^\lambda_r\rt|^p\bigg|\F_s\bigg]^{\frac{1}{p}}\E\bigg[|D_{s,T}|^q
\bigg|\F_s\bigg]^{\frac{1}{q}}\\
\le&\ep M_{s,t}^{\frac{2}{p}}+\frac{1}{\ep}\E\bigg[|D_{s,T}|^q
\bigg|\F_s\bigg]^{\frac{2}{q}}\\
\le&\ep M_{s,t}^{\frac{2}{p}}+\frac{1}{\ep}C_{\delta,q}D_s^2,
\ t_\lambda\le t\le s\le T
\end{aligned}
\end{equation*}
where
\[
M_{s,t}=\E\bigg[\sup_{t\le r\le T}\lt|\overline{Y}_r^\lambda
-Y^\lambda_r\rt|^p\bigg|\F_s\bigg], \ t_\lambda\le t\le s\le T.
\]
Thus, Doob's inequality allows to show that, since $1<p<2$,
\[
\begin{aligned}
&\E\lt[\sup_{s\in[t,T]}M_{s,t}^{\frac{2}{p}}|\F_t\rt]\\
\le&\lt(\frac{2}{2-p}\rt)^{\frac{2}{p}}\E\lt[M_{T,t}^{\frac{2}{p}}|\F_t\rt]\\
\le&\lt(\frac{2}{2-p}\rt)^{\frac{2}{p}}\E\lt[\sup_{s\in[t,T]}
\lt|\overline{Y}^{\lambda}_s-Y^{\lambda}_s\rt|^2|\F_t\rt], \ t\in[t_\lambda,T].
\end{aligned}
\]
Therefore, we can deduce from (\ref{eq:_Y-Y^2}) that
\[
\E\lt[\sup_{s\in[t,T]}\lt|\overline{Y}^{\lambda}_s
-Y^{\lambda}_s\rt|^2|\F_t\rt]\le C_{\delta,\ep} D_t^2
+C_{\delta,\gamma}\ep\lt(\frac{2}{2-p}\rt)^{\frac{2}{p}}
\E\lt[\sup_{s\in[t,T]}\lt|\overline{Y}^{\lambda}_s-Y^{\lambda}_s\rt|^2|\F_t\rt].
\]
By choosing $\ep$ small enough such that
$C_{\delta,\gamma}\ep\lt(\frac{2}{2-p}\rt)^{\frac{2}{p}}<1$,
we get
\[
\E\lt[\sup_{s\in[t,T]}\lt|\overline{Y}^{\lambda}_s
-Y^{\lambda}_s\rt|^2|\F_t\rt]\le C_\delta D_t^2.
\]
Hence, it follows easily from (\ref{eq:_Y-Y^2}) and (\ref{eq:D_t}) that
$$
\begin{aligned}
&\E\lt[\sup_{s\in[t,T]}\lt|\overline{Y}_s^\lambda
-Y^\lambda_s\rt|^2+\int^T_t\lt|\wh{Z}_s^\lambda-Z^\lambda_s\rt|^2\d s
+\int^T_t\int_E \lt|\wh{U}_s^{\lambda}(e)-U_s^\lambda(e)\rt|^2{\mu}(\d
s, \d e)\bigg|\F_t\rt]\\
\le &C_\delta D_t^2.
\end{aligned}
$$
The proof of Lemma \ref{lemma_wtY-Y} is complete now. \hfill$\cajita$
\section*{Acknowledgements}
The author thanks Professor Rainer Buckdahn for his valuable
discussions and advice.

\end{document}